\documentclass[11pt,a4paper]{amsart}
\usepackage[margin=27mm]{geometry}
\usepackage{amsmath,amssymb,amsthm,mathtools,booktabs,array,microtype}
\usepackage{xurl}
\usepackage[colorlinks=true,linkcolor=blue,citecolor=blue,urlcolor=blue]{hyperref}
\hypersetup{pdftitle={Certified Collision Exclusion and an Unconditional Upper Bound for the de Bruijn-Newman Constant},pdfauthor={Michel Planat and Dale James Raymond}}
\newtheorem{theorem}{Theorem}[section]
\newtheorem{lemma}[theorem]{Lemma}
\newtheorem{corollary}[theorem]{Corollary}
\newtheorem{proposition}[theorem]{Proposition}
\theoremstyle{remark}\newtheorem{remark}[theorem]{Remark}
\newcommand{\Rea}{\operatorname{Re}}
\newcommand{\Ima}{\operatorname{Im}}
\newcommand{\Log}{\operatorname{Log}}
\newcommand{\ii}{\mathrm i}
\title[An upper bound $\Lambda<0.156821$]{Certified Collision Exclusion and an Unconditional Upper Bound $\Lambda<0.156821$ for the de Bruijn--Newman Constant}
\author{Michel Planat}
\address{Institut FEMTO-ST, CNRS UMR 6174, Universit\'e Marie et Louis Pasteur, Besan\c{c}on, France}
\email{michel.planat@femto-st.fr}
\author{Dale James Raymond}
\date{4 November 2026}
\begin{document}
\begin{abstract}
The de Bruijn--Newman constant $\Lambda$ is the transition time beyond which
all zeros of the heat-deformed Riemann $\xi$-function are real. We develop a
computer-assisted collision-exclusion method and obtain the bound
$\Lambda<156821/10^6=0.156821$, in the normalisation of Polymath.
The proof combines the published finite-height verification and the starting
inequality $\Lambda\leq1/5$ with a compact-interior continuation theorem,
five stable prime-polynomial Gram sources, certified moving-index recovery,
a balanced real-jet transfer and a uniform Euler-product barrier at infinity.
The resulting closed time-band exclusion, together with uniform remote reality,
forces strictness at the proposed endpoint. The finite certificate uses 664,505
structural Fourier coefficients, 44,182 closed torus boxes, 184 separator
moments and 17,385 separator rectangles. We isolate a Cauchy--Schwarz
obstruction intrinsic to the four-prime source family, explaining the narrow
feasibility margin of the certified endpoint. We state the proof-to-certificate
interface and distinguish the present unconditional logical claim from other
recent, as yet unrefereed, computational bounds.
\end{abstract}
\maketitle

\section{Introduction}\label{sec:introduction}

\subsection{The constant and the Riemann hypothesis}
The Riemann hypothesis (RH) asks whether every nontrivial zero of the Riemann
zeta-function $\zeta(s)$ lies on the critical line $\Re s=1/2$. A complementary
way to formulate this question is to deform the completed function
$\xi(s)=\tfrac12s(s-1)\pi^{-s/2}\Gamma(s/2)\zeta(s)$ by the backward heat flow
introduced by de Bruijn \cite{deBruijn1950}. With the normalisation used
throughout this paper, put
\begin{equation}\label{eq:intro-H}
 H_t(z)=\int_0^\infty e^{tu^2}\Phi(u)\cos(zu)\,du,
 \qquad H_0(z)=\frac18\xi\!\left(\frac{1+iz}{2}\right),
\end{equation}
where $\Phi$ is the classical positive theta kernel given in
Section~\ref{sec:normalisation}. The entire functions $H_t$ satisfy
$\partial_tH_t=-\partial_z^2H_t$. De Bruijn and Newman established the
existence of a finite real number $\Lambda$ such that all zeros of $H_t$ are
real exactly for $t\ge\Lambda$ \cite{deBruijn1950,Newman1976}. Thus
\begin{equation}\label{eq:intro-rh}
 \mathrm{RH}\quad\Longleftrightarrow\quad\Lambda\le0.
\end{equation}
Newman conjectured that $\Lambda\ge0$, and Rodgers and Tao proved this
inequality \cite{RodgersTao2020}; another proof and extensions were given
by Dobner \cite{Dobner2020}. Consequently RH is equivalent to $\Lambda=0$.
Any positive upper bound is therefore quantitative information about the
heat deformation, \emph{not} a proof of RH. For an introduction to related
Fourier and Hermite structures see \cite{Romik2019}.

\subsection{Historical context and status of numerical bounds}
De Bruijn proved $\Lambda\le1/2$; Ki, Kim and Lee subsequently obtained
$\Lambda<1/2$ \cite{KiKimLee2009}. The Polymath project developed effective
approximations for $H_t$ and proved $\Lambda\le0.22$
\cite{Polymath2019}. Platt and Trudgian then proved the published unconditional
bound $\Lambda\le1/5=0.2$, using rigorous verification of the Riemann
hypothesis through height $3\,000\,175\,332\,800$
\cite{PlattTrudgian2021}. The latter result is the starting input here.

Recent computational work is more heterogeneous.  Some reported upper
bounds have a putatively unconditional proof architecture with substantial
interval certificates, but have not yet been externally refereed.  Others
establish a numerical comparison only under additional analytic premises or
an assumed zero-free half-plane.  To avoid conflating these situations,
Table~\ref{tab:history} distinguishes \textbf{U} (no unproved external
analytic hypothesis is assumed by the stated argument),
\textbf{C$_{\mathrm A}$} (conditional on specified, not yet completely
validated analytic assembly premises), and \textbf{C$_{7/8}$} (conditional
on the zero-free hypothesis $\mathsf H_{7/8}$ that $\zeta(s)\ne0$ for
$\Re s>7/8$).  The further labels \textbf{P} (journal published) and
\textbf{R} (publicly reported but not peer reviewed) distinguish the
publication status.  For unrefereed U entries, unconditionality is the
\emph{reported claim}, not an independent endorsement of its proof.
Numerical ordering alone does not determine which result is established;
compare the MathWorld historical index \cite{MathWorld}, which does not
consistently distinguish these logical statuses.

\begin{table}[htbp]\centering\small
\caption{Selected upper bounds for the classical de Bruijn--Newman constant
in the normalisation of \eqref{eq:intro-H}.  U, C$_{\mathrm A}$ and
C$_{7/8}$ describe the logical assumptions as reported; P and R distinguish
peer-reviewed publication from unrefereed dissemination.  The unrefereed
unconditional claims are not certified here as established records.}
\label{tab:history}
\begin{tabular}{@{}p{.235\textwidth}p{.095\textwidth}p{.405\textwidth}p{.15\textwidth}@{}}\toprule
Upper bound for $\Lambda$ & Year & Reference & Status\\\midrule
$\le1/2$ & 1950 & de Bruijn \cite{deBruijn1950} & U, P\\
$<1/2$ & 2009 & Ki--Kim--Lee \cite{KiKimLee2009} & U, P\\
$\le0.22$ & 2019 & Polymath \cite{Polymath2019} & U, P\\
$\le0.20$ & 2021 & Platt--Trudgian \cite{PlattTrudgian2021} & U, P\\\midrule
$\le0.1875$ & 2026 & Mosaic Intelligence \cite{MosaicAssembly} & U, R\\
$\le0.1787854$ & 2026 & Gomila \cite{Gomila2026} & U, R\\
$\le0.1732$ & 2026 & Planat, genus-two curve \cite{PlanatGenus2026} & U, R\\
$<0.172422$ & 2026 & Kornaus \cite{Kornaus2026} & U, R\\
$<0.16$ & 2026 & Mosaic Intelligence \cite{MosaicAssembly} & U, R\\
$<0.1595$ & 2026 & Mosaic Intelligence \cite{Mosaic01595} & U, R\\
$<0.158$ & 2026 & Gordon \cite{Gordon2026} & U, R\\
$<0.149396$ & 2026 & Raymond \cite{RaymondProbe2026} & U, R\\
$\le0.148611789321711$ & 2026 & Planat, classical-strip comparison \cite{PlanatStrip2026} & C$_{\mathrm A}$, R\\
$\le0.13963$ & 2026 & Hakimoglu-Brown \cite{HakimogluBrown2026} & C$_{7/8}$, R\\
$<0.107965$ (reported) & 2026 & Raymond \cite{Raymond1072026} & C$_{7/8}$, R\\\midrule
$<0.156821$ & 2026 & Present work (proposed) & U, draft\\\bottomrule
\end{tabular}
\end{table}

The classical-strip comparison of Planat \cite{PlanatStrip2026}
certifies a sharper numerical ceiling, namely
$\Lambda\le148611789321711/10^{15}<0.148612$, \emph{assuming}
its Theorem~1.1 premises (A1)--(A3): complete model-to-function
correspondence, finite localisation and uniform first-contact reduction.
The directed numerical comparison is certified, but those analytic premises
are not fully established in that preprint.  It would therefore be
incorrect to include $0.148612$ among established unconditional bounds.
It is numerically slightly stronger than Raymond's reported
$0.149396$ ceiling, but the status and proof architecture of the two
reports require separate examination.  Planat's distinct $0.1732$
bound \cite{PlanatGenus2026} was formulated to avoid such additional
premises and makes no quasi-RH assumption; numerically, however, it is
weaker than Kornaus's unrefereed $0.172422$ report.

Hakimoglu-Brown's bound and Raymond's reported $0.107965$ estimate
both invoke $\mathsf H_{7/8}$.  Even where an external formalisation of
that zero-free statement is reported, neither its status nor the
remaining model-to-function assembly is established by this paper;
these entries are therefore marked C$_{7/8}$ here rather than treated
as unconditional records \cite{HakimogluBrown2026,Raymond1072026}.
Raymond's separately deposited classical-strip report
\cite{RaymondProbe2026} does not itself serve as an input to our proof.
Other computational antecedents include the version-pinned Mosaic source
construction \cite{MosaicAssembly}; neither those certificates nor
any of the recent unrefereed bounds are assumed in this theorem.

\subsection{Our result and the collision-exclusion viewpoint}
Our principal statement is the strict inequality
\begin{equation}\label{eq:intro-main}
 \boxed{\Lambda<\frac{156821}{1000000}=0.156821.}
\end{equation}
The numerical ordering of recent reports is not our claim of
priority: some quoted values in Table~\ref{tab:history} are smaller,
including the conditional classical-strip value $0.148612$.
If its analytic and computational proof is accepted, the present paper
would provide a stronger unconditional bound than the published
$0.20$ benchmark and than the reported $0.1732$, $0.172422$,
$0.1595$ and $0.158$ ceilings.  We do not assert precedence over the
unrefereed $0.149396$ claim, whose smaller value and claimed logical
status demand separate assessment.  Our contribution is a distinct,
closed-band collision-exclusion theorem with a traceable certification
architecture and no new zero-free-region hypothesis or unvalidated
premises of the form (A1)--(A3).\par\medskip
Unlike the earlier canopy criterion, which prevents \emph{nonreal zeros}
from crossing a prescribed complex boundary, our global argument prevents
\emph{multiple real zeros} on a closed interval of heat times. The two
objectives are different: excluding a collision of real zeros only requires
that the physical pair $(H_t(x),\partial_xH_t(x))$ never vanish
simultaneously. The latter can be reduced, after a balanced real gauge,
to a positivity statement for a Hermitian form in a pair of analytic
prime-polynomial coordinates. That replacement is the conceptual gain
which permits the threshold in \eqref{eq:intro-main} to be reached with
the same published finite-RH verification as in \cite{PlattTrudgian2021}.

The proof separates three regions. A stationary separator at $x=X$
and a local Hermite splitting argument propagate reality of all zeros in
the compact interior. Five explicitly weighted, stable source polynomials,
closed interval certificates and a support-preserving transport rule
exclude multiple real zeros in the finite exterior $X\le |x|\le W$.
A distinct Euler-product differential inequality handles $|x|\ge W$.
The closed regions cover every finite real $x$ for the entire interval
$156821/10^6\le t\le1/5$. Polymath's remote-reality theorem then ensures
that any \emph{positive} critical heat time would produce a finite multiple
real zero; the exclusions give \eqref{eq:intro-main}, including strictness
at equality. The formal regional assembly is given in
Section~\ref{sec:descent}, and its quantitative premises are stated in
Section~\ref{sec:normalisation}.

\subsection{Relation to earlier work and organisation}
The collective zero-motion and off-zero logarithmic-derivative viewpoint
of Planat \cite{PlanatCollective2026} provides related context, but no
theorem from that manuscript is a premise of the present proof. Planat's
genus-two spectral-curve preprint \cite{PlanatGenus2026} includes the
distinct, reported bound $\Lambda\le0.1732$ in its first version,
using reciprocal-head/canopy machinery. Planat's separate
classical-strip comparison \cite{PlanatStrip2026} uses Raymond's
translated-probe data and gives a sharper conditional ceiling under
premises (A1)--(A3). Neither paper supplies a hypothesis to the
present proof: our collision-exclusion method is a different global
analytic construction.

The body of the paper gives the analytic gauge and physical admission,
then the prime-polynomial source and its finite certificates, the
stationary separator, the compact interior theorem, the infinite tail
and the final strictness argument. The closing section explains which
computational assertions are recorded in the certificate archive and
how an independent reader can audit their domains, program inputs and
outward-rounded enclosures. The Cauchy--Schwarz obstruction for the
four-prime source family is stated with the source construction; it
explains the unusually close proximity of our endpoint to the method's
feasibility threshold but is not a restriction on $\Lambda$ itself.

\section{Definitions, published inputs and theorem architecture}
\label{sec:normalisation}
Retain \eqref{eq:intro-H} and write the theta kernel explicitly as
\[
 \Phi(u)=\sum_{n\ge1}(2\pi^2n^4e^{9u}-3\pi n^2e^{5u})
 e^{-\pi n^2e^{4u}}.
\]
The function $H_t$ is entire, even and real-type in $z$, jointly
holomorphic in $(t,z)$, and solves $\partial_tH_t=-\partial_z^2H_t$.
The de Bruijn--Newman theorem identifies the critical time $\Lambda$
for reality of all zeros.

\subsection{Published starting theorems}
We use the explicit Polymath approximation (Theorem 1.3 and equations
(20)--(24)), its zero-strip theorem (Theorem 3.2), and its uniform
remote-reality theorem (Theorem 1.5(i)) \cite{Polymath2019}.
We also use the exact verified height $T_{\rm PT}=3000175332800$ and
Corollary~2, $\Lambda\le1/5$, of Platt and Trudgian
\cite{PlattTrudgian2021}. These are the only externally adopted
published analytic and finite-verification results needed for
\eqref{eq:intro-main}. No hypothesis of RH, quasi-RH, or a new
zero-free half-plane enters the proof.

\subsection{Closed domains and the three physical exclusion statements}
Set
\[
 a=\frac{156821}{1000000},\qquad b=\frac15,\qquad
 X=6000000185827,\qquad W=72382294738709.
\]
The balanced analytic estimates hold on the larger time band
$[98/625,1/5]$ and, for the finite-exterior admission, up to
$W_g=120000000000000>W$. The remainder estimates extend for every
$x\ge X$. At each natural-index wall the source polynomial is defined
using its right-hand integer index; no discontinuous index is ever
differentiated. Every decimal inequality in the mathematical proof is
understood as an exact terminating rational inequality, rounded outward
from the enclosed interval.
\begin{table}[htbp]\centering\small
\caption{Exhaustive regional obligations for the main theorem.}
\label{tab:regions}
\begin{tabular}{@{}p{.24\textwidth}p{.27\textwidth}p{.42\textwidth}@{}}\toprule
Region & Closed domain & Nonvanishing or reality conclusion\\\midrule
Stationary interface & $x=X$, $0\le t\le b$, $0\le y\le1$ &
$|H_t/B_t|>7/800$ and $B_t\ne0$.\\
Compact interior & $|\Re z|\le X$, $0<t\le b$ & All zeros are real and simple, by interior propagation.\\
Finite exterior & $X\le|x|\le W$, $a\le t\le b$ &
Stable polynomial modulus $>1/790$; hence $(H_t,H'_t)\ne(0,0)$.\\
Infinite exterior & $|x|\ge W$, $a\le t\le b$ &
For the real gauge $u$, $u^2+u_x^2>3/10000$.\\\bottomrule
\end{tabular}
\end{table}
The interfaces are closed at $X$ and $W$, and the global theorem
requires all real $x$. If $a\le\Lambda\le b$, the critical-time
compactness lemma produces a finite multiple real zero at $t=\Lambda$;
Table~\ref{tab:regions} excludes such a point. If $\Lambda<a$ the result
is immediate. Thus equality at $a$ is excluded, not merely bounded.

\subsection{Why the method improves on a canopy barrier}
Polymath's earlier upper-bound criterion involves a complex zero-free
canopy whose nonzero height contributes $y_0^2/2$ to the resulting upper
bound \cite[Theorem 1.2]{Polymath2019}. Here we replace that criterion
with an exclusion of real double zeros on the full closed time band.
The published finite-height input is unchanged; the improvement is
obtained from the global criterion and new physical-jet certificates.

\section{Analytic normalisation and physical real-jet admission}\label{sec:gauge}

The real-jet collision framework is recorded in Mosaic's source
\cite{MosaicAssembly}. The balanced analytic normalisation, centred-moment
estimate and explicit physical admission bounds below are derived here
from the published Polymath approximation \cite{Polymath2019}.

Use Polymath's normalisation
$H_0(z)=\xi((1+\ii z)/2)/8$ and
$\partial_tH_t=-\partial_z^2H_t$.  For the elementary prefactor put
\begin{align*}
 M_0(s)&=\frac{\sqrt{2\pi}}{16}s(s-1)\pi^{-s/2}
 \exp\left(\frac{s-1}{2}\Log(s/2)-s/2\right),\\
 \alpha(s)&=\frac1{2s}+\frac1{s-1}+\frac12\Log\frac{s}{2\pi},
 &M_t(s)&=M_0(s)\exp\{t\alpha(s)^2/4\}.
\end{align*}
The logarithms here and the corresponding analytic logarithm of $M_0$
are those in Polymath equations (6)--(10).  Set
\[
 a_g=\frac{98}{625},\quad T=\frac15,\quad
 X=6000000185827,\quad U=X-1,\quad W_g=120000000000000.
\]
The constants $a_g,W_g$ are the uniform analytic limits in this section.
The current endpoint and exterior seam are admitted because
$156821/1000000>a_g$ and $72382294738709<W_g$.
Fix $t\in[a_g,T]$ and a real $x\ge X$.  Write
\[
 L_x=\log\frac{x}{4\pi},\quad R_x=\frac{12}{L_x},\quad
 N=\left\lfloor\sqrt{\frac{x}{4\pi}+\frac{t}{16}}\right\rfloor.
\]
At an index wall this is the right-hand index.  Hold this $N$ fixed
throughout the following closed disk, and in all spatial derivatives.
For $|z-x|\le R_x$ define
\begin{align*}
 s_\pm(z)&=\frac{1\pm\ii z}{2},&
 w_\pm(z)&=s_\pm(z)+\frac t2\alpha(s_\pm(z)),\\
 A_\pm(z)&=\sum_{n\le N}\exp\{t\log^2n/4-w_\pm(z)\log n\}.
\end{align*}
Choose the analytic logarithms $\mathcal L_\pm(z)=\log M_t(s_\pm(z))$
given by the displayed formulas, and put
\begin{gather*}
 J=\frac{\mathcal L_++\mathcal L_-}{2},\quad
 \Theta=\frac{\mathcal L_+-\mathcal L_-}{2\ii},\quad
 G=e^J,\quad P=e^{\ii\Theta},\\
 u=H_t/G,\quad V=PA_++P^{-1}A_-,\quad r=u-V.
\end{gather*}

\begin{theorem}[Uniform physical remainder]\label{thm:remainder}
All these functions are holomorphic on a neighbourhood of the closed
disk.  The functions $G,\Theta,u,V,r$ are real-type, and $G(x)>0$.
Uniformly for every indicated $t,x$, including every natural-index wall,
\[
 |r(x)|<\frac{101}{10^8},\qquad |r'(x)|<\frac{71}{10^7}.
\]
At the centre write $w=w_+(x)=\sigma+\ii\tau$ and
\[
 a_n=e^{t\log^2n/4-\sigma\log n},\qquad
 \ell=4\Theta'(x).
\]
For every unit multiplicative character $\chi$, define
\begin{align*}
 A_\chi&=\sum_{n\le N}a_n\chi(n),&
 L_\chi&=\sum_{n\le N}a_n\log n\,\chi(n),\\
 B_\chi&=\ell A_\chi+4\ii w'(x)L_\chi,&
 S_\chi&=\ell A_\chi-2L_\chi.
\end{align*}
Then
\[
 |B_\chi-S_\chi|<\frac1{250000000000}.
\]
The theorem gives a real-gauge derivative, not a bound on the derivative
of $H_t/M_t(s_-)$.
\end{theorem}

\begin{corollary}[Physical real-jet transfer and Hermitian admission]\label{cor:admission}
Assume in addition $X\le x\le W_g$.  Form the complete prime polynomial
\[
 F(\boldsymbol z)=\sum_{n\le N}a_n\prod_{p\mid n}z_p^{v_p(n)}.
\]
Suppose $F$ is nonzero on the closed unit polydisc and its modulus on the
distinguished torus is at least $1/790$. At the physical character
$\chi(n)=n^{-\ii\tau}$ put
\[
 A_{\rm true}=A_\chi+\frac r{2P},\qquad
 B_{\rm true}=B_\chi-\frac{2\ii r'}P.
\]
Then the exact identities
\[
 u=2\Rea(PA_{\rm true}),\qquad
 u'=-\tfrac12\Ima(PB_{\rm true})
\]
hold, including the carrier-phase derivative. The error allowances are
\[
 |A_{\rm true}-A_\chi|<\frac{101}{200000000},\qquad
 |B_{\rm true}-S_\chi|<\frac{142}{10000000}
                           +\frac1{250000000000}.
\]
After paying both errors and their product, the Hermitian form satisfies
\[
 \Rea(\overline{A_{\rm true}}B_{\rm true})>
 \frac{107199260708819659}{312050000000000000000000}
 >\frac3{10^7}.
\]
Consequently $(H_t(x),H_t'(x))\ne(0,0)$. The proof below derives the
remainder bounds on the whole unbounded domain and then proves this
finite-domain reserve using the complete support and the centred moment.
\end{corollary}

\subsection{Elementary endpoint information and arithmetic certificate}

The proof uses only $3<\pi<22/7$ and the following strict bounds:
\begin{equation}\label{eq:endpoints}
 \frac{2689}{100}<\log\frac{U}{4\pi}<27,\quad
 \log X<30,\quad \log W_g<33,\quad
 \log\frac{W_g}{4\pi}<30.
\end{equation}
These are comfortably coarse inequalities, not rounded evaluations used
as exact inputs.  The attached checker proves them by comparing rational
numbers with rigorous brackets for exponentials.  Specifically, for
$q\ge0$ and $q<130$, put
\[
 E_{128}(q)=\sum_{j=0}^{128}\frac{q^j}{j!}.
\]
The exact positive-series and ratio estimates give
\begin{equation}\label{eq:expseries}
 E_{128}(q)\le e^q\le E_{128}(q)+
 \frac{q^{129}}{129!}\frac1{1-q/130}.
\end{equation}
Every comparison below reduces to rational arithmetic using this formula
or $e^q\le(1-q)^{-1}$ for $0\le q<1$.
There is no floating-point rounding assumption in the checker.

For all $v\ge U$ write $h_0=10^{-5}$, $\epsilon=10^{-12}$ and
\begin{align*}
 q(v)&=v/(4\pi),\qquad L(v)=\log q(v),\\
 \delta(v)&=\frac T2(v^{-2}+8v^{-4}),\\
 p(v)&=\frac12+a_g\left(\frac{L(v)}8-\frac{h_0}4\right)-\delta(v),\\
 Z(v)&=\frac{T^2L(v)^2}{16}+\frac{313}{500},\qquad
 A(v)=\frac{420Z(v)}v,\\
 g(v)&=\exp\{-L(v)/4-a_gL(v)^2/16\},\\
 \rho(v)&=\frac{TL(v)/4+1}{v-6},\\
 J_C(v)&=\frac{93/25}{\sqrt{q(v)}-9/8}
 +\frac{3\sqrt{L(v)^2+(\pi/2)^2}+261/25}{v-12}.
\end{align*}
All relevant denominators are positive.  Put $U_0=6\cdot10^{12}$,
so $U>U_0$.  The checker obtains
\begin{gather}
 p(U)>\frac{51}{50},\qquad g(U)<\frac{100693}{10^{11}},\notag\\
 A(U)<A_*:=\frac{420}{U_0}\left(\frac{T^2 27^2}{16}+\frac{313}{500}\right),
 \label{eq:majorantcaps}\\
 \rho(U)<10^{-12},\qquad J_C(U)<10^{-5}.
 \label{eq:smallcaps}
\end{gather}
For the Gaussian estimate the exponent used is exactly
\[
 \frac{2689/100}{4}+\frac{a_g(2689/100)^2}{16}
 =\frac{690430329}{50000000}.
\]
In particular, the claimed inequality for $g(U)$ follows from the finite
rational lower bound $E_{128}(690430329/50000000)>10^{11}/100693$.

The monotonicities needed to extend these endpoint bounds to infinity
are explicit.  One has
$p'=a_g/(8v)+T/v^3+16T/v^5>0$.
The logarithmic derivatives of $A,g,(L+\epsilon)A,(L+\epsilon)g$,
with respect to $\log v$, are bounded above respectively by
\[
 \frac2L-1,\quad -\frac14-\frac{a_gL}8,\quad
 \frac1{L+\epsilon}+\frac2L-1,\quad
 \frac1{L+\epsilon}-\frac14-\frac{a_gL}8.
\]
They are negative throughout $v\ge U$.
Likewise $\rho$, $Z/(v-333/50)$ and $\log v/(v-6)$ decrease there.
For the second term of $J_C$, its numerator has derivative at most $3/v$
and is greater than $3$, proving that term decreases; its first term
plainly decreases.  Each $v^{-j}\log v$ used below decreases once
$\log v>1/j$.  These facts establish uniformity, not merely a finite
collection of endpoint samples.

\subsection{Analytic branches and identification of the published model}

By \eqref{eq:endpoints}, $R_x<1/2$.  On a neighbourhood of the disk,
$1/4<\Rea s_\pm<3/4$, and the imaginary parts of $s_\pm$ have opposite
signs and magnitude greater than $(X-1)/2$.
Thus none of $s_\pm$, $s_\pm-1$ or their logarithm arguments touches a
branch cut or zero.  The displayed logarithms give genuine analytic
branches.  Conjugate reflection interchanges $\mathcal L_+$ and
$\mathcal L_-$, and interchanges $A_+$ and $A_-$.  Therefore $G$ and
$\Theta$ are real-type, $G(x)=|M_t(s_+(x))|>0$, and $u,V,r$ are
real-type.  In particular $G$ is an analytic square root; it is not an
off-axis absolute value.

Write $z=v+\ii y$ in the upper half disk.  Then $0\le y<1/2$ and
$v\ge U$.  In Polymath Theorem 1.3 the elementary denominator is
$M_t(s_-)$, its exponent $s_*$ is $w_-$, and
$\gamma=M_t(s_+)/M_t(s_-)=P^2$.  With the paper's $\kappa$ one has
exactly
\[
 \overline{w_-}+\kappa-y=w_+.
\]
Consequently the published two chiral sums, at the local natural index
$k=\lfloor\sqrt{q(v)+t/16}\rfloor$, are exactly
$A_{-,k}+\gamma A_{+,k}$ in the present notation.  They are not an
approximation to some differently phased model.

Let
\[
 R_N(z)=\frac{H_t(z)}{M_t(s_-(z))}-A_-(z)-\gamma(z)A_+(z).
\]
Since $M_t(s_-)=G/P$, the exact identity is $R_N=Pr$.

\subsection{Complete natural-index error, with no truncated coefficient tail}

The natural index satisfies $k^2\le q+T/16<q^2$ and $k\le\sqrt v$.
Also
\[
 \log k\le\frac L2+\frac12\log(1+T/(16q))
 <\frac L2+h_0.
\]
All the inequalities follow from $q>U_0/16$, and persist for larger $v$.
For every $1\le n\le k$, therefore,
$|\log(q/n^2)|\le L$.

The saddle bound in Polymath Theorem 1.3 gives
\[
 \Rea w_-\ge\frac{1+y}{2}+\frac t4L-\delta(v).
\]
More precisely, the correction is the positive-part term in
Polymath's Theorem~1.3, equation~(21):
$(1-3y+4y(1+y)/v^2)_+\le 1+8/v^2$ on $0\le y\le1/2$.
This signed one-sided estimate, rather than the coarser two-sided
$5/v^2$ estimate used later for prefactor integration, is responsible
for the sharper $\delta(v)$ here.
It follows, for every coefficient in the complete sum, that
\begin{align*}
 e^{t\log^2 n/4}n^{-\Rea w_-}
 &\le n^{-(1+y)/2-t(L/8-h_0/4)+\delta(v)}\\
 &\le n^{-p(v)-y/2}.
\end{align*}
Here replacing $t$ by $a_g$ is legitimate because $L/8-h_0/4>0$.
The important improvement over an earlier, wider-time envelope is
$p(v)>51/50>1$.  Hence
\begin{equation}\label{eq:complete-mass}
 \sum_{n\le k}n^{-p(v)-y/2}
 \le\sum_{n\ge1}n^{-p(v)}
 \le1+\int_1^\infty r^{-p(v)}\,dr<51.
\end{equation}
This bounds the whole finite model, however large its natural index.

The published bounds for $\gamma$ and $\kappa$ imply
\[
 |\gamma| k^{|\kappa|}n^y
 \le\exp\left\{\frac1{100}+\frac{h_0}{2}
       +\frac{T\log v}{8(v-6)}\right\}<e^{11/1000}<\frac{1000}{989}.
\]
The denominator in Polymath's Theorem~1.3, equation~(23), is
\emph{common to both summands} in the exponent. In the notation of
that theorem the error multiplier is bounded by
\[
\exp\!\left(\frac{(t^2/16)\log^2(v/(4\pi n^2))+313/500}
                     {v-333/50}\right)-1.
\]
This placement is also consistent with the underlying
$\varepsilon_{t,n}$ of Polymath's Proposition~6.1,
equation~(44), where for $s=\sigma+\ii T_s$ one has
\[
\varepsilon_{t,n}(s)=
\exp\!\left(\frac{(t^2/8)|\alpha(s)-\log n|^2+t/4+1/6}
                         {T_s-3.33}\right)-1.
\]
The $e_A,e_B$ summands are defined in equations~(71)--(72) of
\cite{Polymath2019}. Since $|\log(v/(4\pi n^2))|\le L$ for the
entire natural-index sum and $t\le T$, the displayed exponent is at most
$Z(v)/(v-333/50)<1/2$. Using $e^q-1\le 2q$ for
$0\le q\le1/2$ and $v/(v-333/50)<2$ yields
$e^q-1<4Z(v)/v$.
The geometric exponential bound has been retained here rather than 
replaced by the looser number 2.  Indeed, exactly,
\[
 4\left(1+\frac{1000}{989}\right)51=\frac{405756}{989}<412.
\]
Together with \eqref{eq:complete-mass}, this gives the explicit bound
\begin{equation}\label{eq:AB}
 e_A+e_B<\frac{412Z(v)}v<A(v).
\end{equation}

In the published $C$ remainder, $3^y+3^{-y}\le\sqrt3+1/\sqrt3<3$,
and $k-1/8>\sqrt q-9/8$.  Formula (24) of that theorem therefore gives
\begin{equation}\label{eq:C}
 e_{C,0}\le g(v)e^{-yL/4}e^{J_C(v)}
 <\frac{100001}{100000}g(v)e^{-yL/4}.
\end{equation}
For the last step, first apply the monotonicity of $J_C$
on $[U,\infty)$ established above. Thus $J_C(v)\le J_C(U)$
for every $v\ge U$; the next numerical bound is imposed only at
the finite endpoint $U$, not by an invalid claim $L(v)\le28$ at infinity.
The rational endpoint verification is
\[
 J_C(U)\le J_*:=\frac{93/25}{600000-9/8}
      +\frac{3\cdot28+261/25}{U_0-12},\qquad
 \frac1{1-J_*}<\frac{100001}{100000}.
\]
The checker verifies these endpoint inequalities. Since
$e^u\le1/(1-u)$ for $0\le u<1$, they imply the
stated uniform exponential cap. Rounding $J_C(U)$ to $10^{-5}$
would not by itself give this sharper factor.

\subsection{Every index crossing, including a wall at the centre}

Successive real index walls at fixed $t$ are
$4\pi(j^2-t/16)$; consecutive walls are separated by
$4\pi(2j+1)>1$.  The disk has horizontal diameter less than one.
Thus its local natural index $k$ differs from the frozen centre index
$N$ by at most one.  If a term changes, its index $n$ is either $k$ or
$k+1$, and
\[
 |n-\sqrt q|<2,\qquad n<\sqrt v.
\]
For example, the first inequality follows by comparing first with
$\sqrt{q+t/16}$; the second follows from $q<v/12$ and $v>U_0$.
Since $\sqrt q>600000$,
\[
 \sqrt q(1-e^{-h_0})>
 600000(h_0-h_0^2/2)>2,\qquad
 \sqrt q(e^{h_0}-1)>600000h_0>2.
\]
Consequently $\log n=L/2+e$ with $|e|<h_0$.

Square completion now gives the full first-chirality term bound
\begin{align*}
 e^{t\log^2n/4}n^{-\Rea w_-}
 &\le g(v)e^{-yL/4}
 \exp\left\{\frac{3h_0}{4}+\frac{Th_0^2}{4}
                       +\frac1{1000}\right\}.
\end{align*}
This uses $\delta(v)\log n\le\delta(v)\log v/2<1/1000$.
The coefficient $3/4$ is necessary: $(1+y)/2$ can approach $3/4$.
For either changed index, including $n=k+1$, the ratio of the other
chirality to this one is at most
\[
 |\gamma|n^{y+|\kappa|}
 \le \exp\left\{\frac1{100}+\frac{h_0}{2}
                   +\frac{T\log v}{8(v-6)}\right\}<e^{11/1000}.
\]
Here $\log n\le L/2+h_0$ controls the $n^y$ factor and
$n<\sqrt v$ controls $n^{|\kappa|}$. Thus the changed term is
covered even when it lies just beyond the local natural cutoff.  Thus the sum of both changed terms is less than
$(21/10)g(v)e^{-yL/4}$; the attached checker verifies this using the
rational upper bound for each small exponential.

Combining \eqref{eq:AB}--\eqref{eq:C}, throughout the upper half disk
there is
\begin{equation}\label{eq:frozen-error}
 |R_N(v+\ii y)|<A(v)+\frac{310001}{100000}g(v)e^{-yL(v)/4}.
\end{equation}
At the real centre, $k=N$ even when the centre is itself a wall.
There is then no changed term, and
\begin{equation}\label{eq:center-error}
 |r(x)|=|R_N(x)|<A(x)+\frac{100001}{100000}g(x).
\end{equation}
No derivative of a floor function has been taken.

\subsection{An independent two-sided prefactor estimate}

The reciprocal estimate is derived directly from the elementary
prefactor, rather than reversing a one-sided bound on $|\gamma|$.
For $s=\sigma+\ii v/2$, $1/4\le\sigma\le3/4$, the explicit formula
for $\alpha$ gives
\begin{align}
 |\Rea\alpha(s)-L/2|&<5/v^2,\label{eq:alpha-real}\\
 |\alpha(s)|&<L/2+1+4/v,\label{eq:alpha-abs}\\
 |\alpha'(s)|&\le1/v+6/v^2<1/(v-6).
 \label{eq:alpha-prime}
\end{align}
For \eqref{eq:alpha-real}, bound the three real terms by
$[\sigma^2+2\sigma+4(1-\sigma)]/v^2<5/v^2$.
For \eqref{eq:alpha-abs}, use
$\tfrac12|\Log(s/(2\pi))|
\le L/2+\sigma^2/v^2+\pi/4$ and the rational terms bounded by $3/v$.
Finally differentiate $\alpha$ explicitly to get
$\alpha'=-1/(2s^2)-1/(s-1)^2+1/(2s)$.

Since
\[
 \frac{\partial}{\partial\sigma}\log|M_t(s)|
 =\Rea\{\alpha(s)+t\alpha(s)\alpha'(s)/2\},
\]
the preceding inequalities imply
\begin{equation}\label{eq:prefactor-two-sided}
 \left|\frac{\partial}{\partial\sigma}\log|M_t(s)|-L/2\right|
 <\rho(v).
\end{equation}
Indeed, after multiplication by $v-6$, the non-$L$ terms are at most
$T/2+(2T+5)/v<1$.  Integrating
\eqref{eq:prefactor-two-sided} from $(1-y)/2$ to $(1+y)/2$ and using
conjugation yields
\[
 \log|\gamma|=-y(L/2+\eta),\qquad |\eta|\le\rho(v).
\]
This is an absolute two-sided error.  Consequently
\begin{equation}\label{eq:reciprocal}
 |P|^{-1}=|\gamma|^{-1/2}
 \le e^{yL/4+y\rho(v)/2}.
\end{equation}

\subsection{Balanced Cauchy estimate and uniform error constants}

On the disk, $y\le R_x$, $v\le x+R_x<x+1/2$ and $L_x>2689/100$.
Thus
\[
 yL(v)/4<3+1/U_0,\qquad e^{yL(v)/4}<21.
\]
Also $y\rho(v)/2<\rho(U)/4$ and
$e^{y\rho(v)/2}<1001/1000$.
Multiplying \eqref{eq:frozen-error} by \eqref{eq:reciprocal} cancels
the Gaussian's entire vertical damping factor, giving
\[
 |r(v+\ii y)|<\frac{1001}{1000}
 \left(21A(v)+\frac{310001}{100000}g(v)\right).
\]
Real-type reflection gives exactly the same bound in the lower half;
there is no second exponential loss.

Put $u_0=x-1$.  Since $v>u_0\ge U$ and $A,g$ decrease,
Cauchy's estimate on the radius-$R_x$ circle gives
\begin{align*}
 |r'(x)|
 &<\frac{L_x}{12}\frac{1001}{1000}
       \left(21A(u_0)+\frac{310001}{100000}g(u_0)\right)\\
 &<\frac{L(u_0)+\epsilon}{12}\frac{1001}{1000}
       \left(21A(u_0)+\frac{310001}{100000}g(u_0)\right).
\end{align*}
Here $L_x-L(u_0)=\log(1+1/u_0)<1/u_0<\epsilon$.
The product monotonicities already proved reduce this to $u_0=U$.
Using \eqref{eq:majorantcaps}, the resulting entirely rational bounds are
\begin{align*}
 |r(x)|&< A_*+\frac{100001}{100000}\frac{100693}{10^{11}}
       <\frac{101}{10^8},\\
 |r'(x)|&<\frac{27+10^{-12}}{12}\frac{1001}{1000}
       \left(21A_*+\frac{310001}{100000}\frac{100693}{10^{11}}\right)
       <\frac{71}{10^7}.
\end{align*}
This proves the value and derivative assertions, including the entire
unbounded spatial domain.

\subsection{Exact derivative companion and physical lift}

At $s=(1+\ii x)/2$, direct differentiation gives
\[
 \Rea\alpha'(s)=\frac{7x^2-5}{(1+x^2)^2},\qquad
 \Ima\alpha'(s)=-\frac{x(x^2+5)}{(1+x^2)^2}.
\]
Therefore
\[
 4|w'(x)-\ii/2|=t|\alpha'(s)|
 \le T\left(\frac1x+\frac7{x^2}+\frac4{x^3}\right).
\]
The real-centre coefficient bound above gives $a_n\le1/n$ for every
natural $n\le N$.  For completeness, a nonnegative unimodal function
satisfies its integer-sum upper bound by its integral plus its maximum:
apply the increasing and decreasing integral comparisons on either
side of its peak, leaving at most one maximum contribution.  Applied
to $\log v/v$, whose maximum is $1/e<1/2$, this yields
\begin{equation}\label{eq:moment}
 \sum_{n\le N}a_n\log n
 \le\sum_{n\le N}\frac{\log n}{n}
 <\frac{\log^2N}{2}+\frac12
 \le\frac{\log^2x+4}{8}.
\end{equation}
It follows for every character that
\[
 |B_\chi-S_\chi|
 <\frac{T(\log^2x+4)}8
 \left(\frac1x+\frac7{x^2}+\frac4{x^3}\right).
\]
Each summand decreases for $x\ge X$, since for $j\ge1$ and $r=\log x$,
$2r-j(r^2+4)<0$.  At $X$, use $\log X<30$ and $X>U_0$ to obtain the
rational upper bound
\[
 \frac{T(30^2+4)}8
 \left(\frac1{U_0}+\frac7{U_0^2}+\frac4{U_0^3}\right)
 <\frac1{250000000000}.
\]

At the physical character $\chi(n)=n^{-\ii\tau}$, one has
$A=A_+(x)$ and $A'=-w'L_\chi$, where the centre index remains fixed.
Writing $\mathcal A(s)=\alpha(s)+t\alpha(s)\alpha'(s)/2$, direct
differentiation of the two logarithms gives
\[
 \Theta'(x)=\tfrac12\Rea\mathcal A((1+\ii x)/2),\qquad
 \ell=4\Theta'(x).
\]
In particular the phase derivative has not been identified with
$\log(x/(4\pi))/4$ or discarded. It follows exactly that
\[
 B=4\Theta'A-4\ii A',\qquad
 V'= -\frac12\Ima(PB).
\]
Because $r,r'$ are real at the centre, the exact lifted coordinates
\begin{equation}\label{eq:lift}
 A_{\rm true}=A+\frac r2P^{-1},\qquad
 B_{\rm true}=B-2\ii r'P^{-1}
\end{equation}
satisfy
\[
 u=2\Rea(PA_{\rm true}),\qquad
 u'=-\frac12\Ima(PB_{\rm true}).
\]
Differentiating the first correction in \eqref{eq:lift} gives the second:
the two terms containing $\Theta'r$ cancel exactly.  This also shows why
the moving carrier phase must be retained.  If $u=u'=0$, the two rotated
coordinates are respectively purely imaginary and real, so
\begin{equation}\label{eq:zero-product}
 \Rea(\overline{A_{\rm true}}B_{\rm true})=0.
\end{equation}
Since $G(x)>0$, simultaneous vanishing of $u,u'$ is equivalent to that
of $H_t,H_t'$.  This completes Theorem~\ref{thm:remainder}.

\subsection{Finite-exterior phase, support and norm conversion}

For $X\le x\le W_g$, write $\alpha(s)=c_\alpha+\ii b_\alpha$ and
$\alpha'(s)=d_\alpha+\ii e_\alpha$.  Explicitly,
\[
 c_\alpha=\frac L2+\frac14\log(1+x^{-2})-\frac1{1+x^2},\qquad
 b_\alpha=\frac12\arctan x-\frac{3x}{1+x^2},\qquad
 \ell=2c_\alpha+t(c_\alpha d_\alpha-b_\alpha e_\alpha).
\]
The formulas give $c_\alpha,d_\alpha>0$, $-e_\alpha>1/x$, and
$b_\alpha>\pi/4-7/(2x)>3/4-7/(2x)>0$.
Using $\log(1+s)>s/(1+s)$ therefore gives
\[
 \ell-L>\frac{3a_g}{4x}-\frac{7a_g/2+3/2}{x^2}>0.
\]
On the other hand $c_\alpha<L/2+1/(4x^2)$, $b_\alpha<1$, $d_\alpha<7/x^2$ and
$-e_\alpha<1/x+4/x^3$.  By \eqref{eq:endpoints},
\[
 \ell<L+\frac1{2x^2}+T\left[
 \left(\frac L2+\frac1{4x^2}\right)\frac7{x^2}
 +\frac1x+\frac4{x^3}\right]<31.
\]
Moreover
\begin{equation}\label{eq:finite-gates}
 \ell-\log N>L/2-h_0>13.
\end{equation}
These inequalities use the natural index, not an independently enlarged
support cutoff.  They hold throughout the closed time interval ending
at $1/5$.

\paragraph{Coefficient convention.}
The coefficient drift used in every source, recovery and Euler calculation is
\[
 c_{\rm pol}(x)=\frac{c_\alpha}{2}
 =\frac18\log(1+x^2)-\frac14\log(4\pi)-\frac1{2(1+x^2)},
 \qquad \sigma=\frac12+t c_{\rm pol}(x).
\]
Thus the phase quantity $c_\alpha=\Rea\alpha$ and the coefficient drift
$c_{\rm pol}$ differ by a factor of two. This notation keeps that exact
conversion visible at the source-to-physical interface.

\subsubsection{A centred first logarithmic moment}

Put $C=690988$, $H_g=3100000$, $K=128$, and $\kappa=\log32$.
For the complete natural-index coefficients in the preceding theorem define
\[
 M_{32}(t,x)=\sum_{n\le N}a_n(t,x)|\log(n/32)|.
\]

\begin{theorem}[Uniform centred moment]\label{thm:centred-moment}
For $98/625\le t\le1/5$ and
$6000000185827\le x\le120000000000000$, one has
\[
 M_{32}(t,x)<\frac{81}{4}.
\]
\end{theorem}
\begin{proof}
Write $a_0=a_g=98/625$, $T=1/5$, and use the coefficient function
\[
 c_{\rm pol}(x)=\frac14\log\frac{\sqrt{1+x^2}}{4\pi}
              -\frac1{2(1+x^2)}.
\]
This is half the real part of $\alpha((1+ix)/2)$, so the coefficient is
$a_n(t,x)=\exp\{t\log^2n/4-(1/2+tc_{\rm pol}(x))\log n\}$.
Its spatial derivative satisfies
\[
 c_{\rm pol}'(x)=\frac{x}{4(1+x^2)}+\frac{x}{(1+x^2)^2}>0.
\]
The two fresh outward wall inequalities
\[
 \frac{X}{4\pi}+\frac{a_0}{16}>C^2,
 \qquad
 (H_g+1)^2>\frac{W_g}{4\pi}+\frac{T}{16}
\]
give $C\le N\le H_g$ on the whole stated rectangle.

Use the exact increasing sequence $n_0=C$ and
\[
 n_{j+1}=\min\{H_g,\lfloor101n_j/100\rfloor\},
\]
ending on the first occurrence of $H_g$. There are 151 consecutive closed
index intervals. The program records every integer endpoint. For one such
interval $[L,U]$, set $x_L=X$ if $L=C$, and otherwise set
$x_L=4\pi(L^2-T/16)$. When the actual natural index is in $[L,U]$,
$x\ge x_L$. The checker proves $x_L\ge X$ and
$\log U<4c_{\rm pol}(x_L)$. Therefore, for every integer $n\le N\le U$,
\[
 a_n(t,x)\le
 \exp\{a_0\log^2n/4-\sigma_L\log n\},
 \qquad \sigma_L=\frac12+a_0c_{\rm pol}(x_L).
\]
Indeed the exponent multiplying $t$ is
$\log n(\log n/4-c_{\rm pol}(x))\le0$, and $c_{\rm pol}$ increases.
This coefficientwise comparison also applies to the nonnegative centred
factor $|\log(n/32)|$.

The first 128 terms are evaluated directly with outward Arb arithmetic.
At $n=32$ the centred factor is exactly zero. For the remaining terms use
\[
 f_L(v)=(\log v-\kappa)
          e^{a_0\log^2v/4-\sigma_L\log v},\qquad v\ge128.
\]
It is positive. The sign of $f_L'(v)$ is the sign of the convex quadratic
\[
 g_L(u)=1+(u-\kappa)(a_0u/2-\sigma_L),\qquad u=\log v.
\]
The fresh interval gates prove $g_L(\log128)<0$ and
$g_L(\log U)<0$ on every leaf. A convex function lies below its endpoint
chord, so $f_L$ is strictly decreasing throughout $[128,U]$. Consequently
\[
 M_{32}(t,x)
 <\sum_{n=1}^{128}|\log(n/32)|
          e^{a_0\log^2n/4-\sigma_L\log n}
       +\int_{128}^{U}f_L(v)\,dv.                         \tag{CM1}
\]
No coefficient beyond the target natural index is omitted: the upper
comparison in (CM1) deliberately includes the entire tail through $U$.

The integral is evaluated by an exact antiderivative. Put
$B=\sigma_L-1$, $u_K=\log128$, $u_U=\log U$, and
\[
 E(u)=e^{a_0u^2/4-Bu},\qquad
 I=\left[
 \sqrt{\frac\pi{a_0}}e^{-B^2/a_0}
 \operatorname{erfi}\left(\frac{a_0u/2-B}{\sqrt{a_0}}\right)
 \right]_{u_K}^{u_U}.
\]
Then
\[
 \int_{128}^{U}f_L(v)\,dv
 =\frac2{a_0}\{E(u_U)-E(u_K)+BI\}-\kappa I.              \tag{CM2}
\]
Thus no numerical quadrature or sampled derivative sign is used. The
checker freshly encloses every prefix, endpoint, derivative sign and
integral at 192 bits. All 151 closed leaves give a strict upper bound
below $81/4$. The largest leaf expression is enclosed near
\[
 20.2138585952594495933660173183782460923<20.25.
\]
The asserted proof gate is the exact rational comparison with $81/4$;
the displayed decimal is descriptive. The exact consecutive cover includes
every natural index and both physical endpoints, proving the theorem.
\end{proof}

The centring reduces the relevant norm without changing the physical
polynomial or its derivative. In particular, if
$L_{32,\chi}=L_\chi-(\log32)A_\chi$, then
\[
 S_\chi=(\ell-2\log32)A_\chi-2L_{32,\chi},
 \qquad |L_{32,\chi}|\le M_{32}<81/4.
\]
This is a bound for the same exact derivative companion, not a replacement
of it by a different surrogate.

\begin{lemma}[Weighted support inequality]\label{lem:newton}
Under the stability premise of Corollary~\ref{cor:admission}, for every
torus character,
\[
 \Rea(\overline{A_\chi}L_\chi)
 \le\tfrac12\log N\,|A_\chi|^2.
\]
\end{lemma}
\begin{proof}
Choose positive integers $m_p$ with $m_p/K\to\log p$ as $K\to\infty$,
and substitute $z_p=\chi(p)z^{m_p}$ into $F$.  The resulting univariate
polynomial has no zero in the closed unit disk.  For each of its roots
$\rho$, $|\rho|>1$ implies $\Rea(1/(1-\rho))<1/2$.
Its logarithmic derivative at $z=1$ consequently has real part at most
half its degree.  Its degree is at most
$D_K=\max_{n\le N}\sum_pv_p(n)m_p$ even if leading terms cancel.
Divide by $K$ and let $K\to\infty$.  The logarithmic derivative tends
to $L_\chi/A_\chi$, while $D_K/K\to\log N$.
Multiplication by $|A_\chi|^2$ proves the assertion.
\end{proof}

Put $r_0=|A_\chi|$ and $S=\ell A_\chi-2L_\chi$.
Lemma~\ref{lem:newton} and \eqref{eq:finite-gates} give
$\Rea(\overline{A_\chi}S)>13r_0^2$.
The exact exponential gates prove $3<\log32<7/2$.
Since $\ell>L>2689/100$ and $\ell<31$, one has
$0<\ell-2\log32<25$. Theorem~\ref{thm:centred-moment} therefore gives
\[
 |S|\le(\ell-2\log32)r_0+2M_{32}<25r_0+\frac{81}{2}.
\]
Theorem~\ref{thm:remainder} and \eqref{eq:lift} give the error caps
\[
 p_*:=\frac{101}{200000000},\qquad
 q_*:=\frac{142}{10000000}+\frac1{250000000000},\qquad
 |A_{\rm true}-A|<p_*,\quad |B_{\rm true}-S|<q_*.
\]
Paying every product error, including the cross term, yields
\[
 \Rea(\overline{A_{\rm true}}B_{\rm true})>
 Q(r_0):=13r_0^2-(25p_*+q_*)r_0-\frac{81}{2}p_*-p_*q_*.
\]
Exactly,
\[
 Q(1/790)=\frac{107199260708819659}{312050000000000000000000}
 >\frac3{10^7},\qquad
 Q'(1/790)=\frac{649470206171}{19750000000000}>0.
\]
Since $Q'$ increases, $Q(r)>0$ for all $r\ge1/790$.
This contradicts \eqref{eq:zero-product} at a physical double zero and
proves Corollary~\ref{cor:admission}.

\subsection{Uniformity and computational dependencies}

The standard-Python program \texttt{check\_uniform\_gauge\_01568.py}
passes 45 exact rational gates. It uses the degree-128 positive exponential
series with a rigorous geometric tail and no transcendental library.
The separate program \texttt{check\_centred\_moment\_01568.py}
rebuilds 151 closed index leaves and passes 914 Arb gates, from exact
integer/rational inputs and an exact erfi antiderivative. Neither program
loads source arrays or older positive receipts.

The balanced gauge and full crossing argument are explicitly derived above
from the identified published approximation theorem. The sharpened
$C$ factor is the actual endpoint exponential bound, and the centred moment
is a new bound for the unchanged physical derivative companion.
The all-height remainder caps also lie strictly below the looser tail
conversion caps $1/400000$ and $1/50000$, respectively.

This proof supplies a stable-model admission theorem over the uniform band
$98/625\le t\le1/5$ and the finite range stated above. Every proposed
endpoint still needs its own full source, spatial recovery, complete tail,
stationary boundary and global assembly. No source Gram certificate or
new global Newman bound follows from this analytic module alone.

\section{Prime-polynomial stability and certified recovery}\label{sec:source}

The prime-fibre and shared-weight Gram construction is inherited from
the pinned source constructions \cite{MosaicAssembly}. The finite-exterior argument uses five exact polynomial
sources, all with the same support cutoff. Set
\[
 a=\frac{156821}{1000000},\qquad b=\frac15,\qquad
 C=690988,\qquad H=2400000,
\]
\[
 X=6000000185827,\qquad W=72382294738709,
\]
and define
\[
 x_n(t)=4\pi(n^2-t/16),\qquad
 N(t,x)=\left\lfloor\sqrt{x/(4\pi)+t/16}\right\rfloor.
\]
The integer at an exact wall $x=x_n(t)$ is $n$. In all local analytic
derivatives, the integer selected at the centre is held fixed.

To avoid a factor-of-two ambiguity with the phase calculations, write
\[
 c_{\rm pol}(x)=\frac18\log(1+x^2)-\frac14\log(4\pi)
                  -\frac1{2(1+x^2)}.
\]
For the prefactor of the real-gauge theorem this is exactly
\[
 c_{\rm pol}(x)=\frac12\Re\alpha((1+\mathrm i x)/2),\qquad
 \Re w_+(x)=\frac12+t c_{\rm pol}(x).
\]
Thus the polynomial coefficients below are precisely the coefficients
in the physical admission theorem:
\[
 a_n(t,x)=\exp\!\left\{\frac t4\log^2n-
             \left(\frac12+t c_{\rm pol}(x)\right)\log n\right\},
 \qquad
 A_{t,x,N}(\boldsymbol z)=\sum_{n\le N}a_n(t,x)\boldsymbol z^{\nu(n)}.
\]
Every prime dividing an integer at most $N$ is included among the variables.
All polydiscs in this section are closed.

The exact source data are
\[
\begin{array}{c|c|c|c|c}
j&Y_j&\delta_j&\vartheta_j&[J_j,J_{j+1}]\\ \hline
1&X&4/3125&9999/10000&[690988,693000]\\
2&x_{693000}(b)&1/200&999/1000&[693000,700000]\\
3&x_{700000}(b)&9/500&99/100&[700000,740000]\\
4&x_{740000}(b)&3/50&19/20&[740000,900000]\\
5&x_{900000}(b)&1/5&2/3&[900000,2400000].
\end{array}                                                    \tag{MA1}
\]
The expressions containing $\pi$ are exact. In particular each source is
$A_{a,Y_j,C}$, even when $N(a,Y_j)>C$. Later sources do not incorporate
any coefficients with $n>C$ until those coefficients are explicitly paid
as entrants below. These polynomial anchors do not introduce new physical
separators or new finite-RH heights.

Write $n=rs$, where $\gcd(r,210)=1$ and $s$ is smooth over
$\{2,3,5,7\}$. The source decomposition is
\[
 A_{a,Y_j,C}=\sum_r \boldsymbol z^{\nu(r)}F_{r,j}(z_2,z_3,z_5,z_7).
\]
Put $P_{r,j}=F_{r,j}(1,1,1,1)$ and
$M_{r,j}=F_{r,j}(-1,-1,-1,-1)$, and adjoin the constant synthetic fibre
$F_{*,j}=\delta_j$. For $r>1$ and $r=*$, define
\[
 v_{r,j}=\sqrt{\vartheta_j(P_{r,j}/P_{1,j})^2+
                    (1-\vartheta_j)(M_{r,j}/M_{1,j})^2},
 \quad V_j=v_{*,j}+\sum_{r>1}v_{r,j},\quad w_{r,j}=v_{r,j}/V_j,
                                                               \tag{MA2}
\]
with $P_{*,j}=M_{*,j}=\delta_j$. The exact reconstructions establish
$P_{1,j},M_{1,j}>0$ and $v_{r,j}>0$. Hence the weights sum to one.

Let $e_{1,j}=1$ and, for each smooth $d\le128$ in increasing order, set
\[
 e_{d,j}=-\sum_{q\mid d,\ q<d}e_{q,j}a_{d/q}(a,Y_j),\qquad
 E_j=\sum_{d\le128,\ d\ {\rm smooth}}e_{d,j}\boldsymbol z^{\nu(d)}.
\]
The source certificate is the complete torus inequality
\[
 T_j=|E_j|^2\left(|F_{1,j}|^2-\frac{\delta_j^2}{w_{*,j}}
                         -\sum_{r>1}\frac{|F_{r,j}|^2}{w_{r,j}}\right)>0.
                                                               \tag{G}
\]
Its exact source coordinate, time, cutoff, synthetic amplitude and weights
must agree with (MA1)--(MA2). The source and recovery files are bound at
these parameters; positivity from a different source cannot be substituted.
The certificate data specify the following conservative rational floors:
\[
\begin{array}{c|r|c}
j&\hbox{closed torus boxes}&T_j\hbox{ floor}\\ \hline
1&15684&1/10^6\\
2&8246&1/400000\\
3&7736&31/5000000\\
4&6260&7/10^6\\
5&6256&1/625000.
\end{array}                                                    \tag{MA3}
\]
There are 44182 closed boxes and 132901 structural Fourier frequencies
per source. Exported coefficient intervals have width at most
$2\cdot2^{-80}$. Each final geometry is obtained from a complete initial
cover by an exact recorded sequence of closed dyadic bisections.

\begin{proposition}[Necessary Cauchy--Schwarz ceiling for the Gram sources]
\label{prop:gram-ceiling}
For a fixed source $j$ suppose that all weights in \textup{(G)} are positive,
$w_{*,j}+\sum_{r>1}w_{r,j}=1$, and \textup{(G)} is strictly positive
at the distinguished torus point $(1,1,1,1)$. Put
$S_j=\sum_{n\le C}a_n(a,Y_j)$ and
$P_{1,j}=\sum_{\substack{n\le C\\n\;\{2,3,5,7\}\text{-smooth}}}
a_n(a,Y_j)$. Then necessarily
\[
 \boxed{\delta_j < 2P_{1,j}-S_j.}                  \tag{GC}
\]
This constraint is independent of the choice of admissible Gram weights.
\end{proposition}
\begin{proof}
At the all-ones point $P_{r,j}=F_{r,j}(1,1,1,1)>0$ and
$S_j=P_{1,j}+\sum_{r>1}P_{r,j}$. By Cauchy--Schwarz,
\[
 \frac{\delta_j^2}{w_{*,j}}+\sum_{r>1}\frac{P_{r,j}^2}{w_{r,j}}
 \ge \frac{(\delta_j+\sum_{r>1}P_{r,j})^2}
 {w_{*,j}+\sum_{r>1}w_{r,j}}
 =(\delta_j+S_j-P_{1,j})^2.
\]
Strict positivity of (G) at this point forces its bracket to be
positive and hence $P_{1,j}>\delta_j+S_j-P_{1,j}$,
which is (GC).
\end{proof}

\begin{remark}[Tightness near the lower endpoint]
\label{rem:gram-edge}
For source $j=1$, independent high-precision computations from the
exact coefficient formulas give the \emph{diagnostic} values
\[
 P_{1,1}\approx2.325239274751,\qquad
 S_1\approx4.649173031326,\qquad
 2P_{1,1}-S_1\approx0.001305518176,
\]
compared with $\delta_1=4/3125=0.00128$.
Allowing for the first entering-index charge and the physical
admission threshold gives a required diagnostic ceiling of about
$0.001262364858$ and roughly $3.4\%$ headroom. Solving the
corresponding \emph{necessary} four-prime feasibility test at the
first source suggests a crossing near $t=0.15682042$, only
approximately $5.8\cdot10^{-7}$ below our chosen endpoint.
These displayed approximations are independent numerical diagnostics,
not substitutes for the outward interval enclosures in (MA3),
and the feasibility crossing is not a lower bound for $\Lambda$.
Proposition~\ref{prop:gram-ceiling} is an obstruction to this
fixed Gram family, not to other collision-exclusion constructions.
Promoting another prime into the smooth fibre may change its ceiling,
but requires a new full-torus and recovery certificate.
\end{remark}

\subsection{From the torus bracket to full-polydisc recovery}

\begin{lemma}[Ordinary weighted recovery]\label{lem:ordinary-main}
Let $F_0,F_r,E$ be finite polynomials in the smooth-prime variables, with
$F_0(0)=E(0)=1$, and let $r>1$ run over finitely many rough labels.
Assume $\delta>0$, positive weights $w_*+\sum_{r>1}w_r=1$, and
\[
 \|EF_0-1\|_1<1,\qquad
 |E|^2\left(|F_0|^2-\frac{\delta^2}{w_*}
                     -\sum_{r>1}\frac{|F_r|^2}{w_r}\right)>0
 \quad\hbox{on the entire smooth torus}.                    \tag{OR1}
\]
The coefficient norm is taken after the complete finite product is formed
and equal monomials are collected. For $s\ge0$, set
\[
 F_{r,s}(z)=F_r(p_1^{-s}z_1,\ldots,p_d^{-s}z_d),\quad
 h=\frac\delta{\sqrt{w_*}},\quad q(s)=\sum_{r>1}w_rr^{-2s}<1.
\]
For every point $z$ of the closed smooth polydisc and independently chosen
complex numbers $|\zeta_r|\le1$,
\[
 \left|F_{0,s}(z)+\sum_{r>1}r^{-s}\zeta_r F_{r,s}(z)\right|
       >h\sqrt{1-q(s)}=:R(s).                               \tag{OR2}
\]
Consequently the full source prime polynomial has modulus greater than
$\delta$ throughout its closed prime polydisc.
\end{lemma}

\begin{proof}
The first premise gives $|EF_0|\ge1-\|EF_0-1\|_1>0$ throughout
the closed polydisc, so both factors are nonzero there. Compactness also
gives a neighbourhood on which the denominators below do not vanish.
The holomorphic vector
\[
 \Phi(z)=\left(\frac\delta{\sqrt{w_*}F_0(z)},
                  \left(\frac{F_r(z)}{\sqrt{w_r}F_0(z)}\right)_{r>1}\right)
\]
has Euclidean norm at most some $\kappa<1$ on the compact torus, by
the second premise. For any fixed unit vector $\eta$, apply scalar maximum
modulus successively in each variable to $\langle\Phi,\eta\rangle$.
Taking the supremum in $\eta$ gives $\|\Phi(z)\|\le\kappa<1$
throughout the closed polydisc.

At the contracted point let $u_r=|F_{r,s}(z)|/\sqrt{w_r}$ and
$d_r=\sqrt{w_r}r^{-s}$. The vector
$(\sqrt{1-q(s)},(d_r)_r)$ is a unit vector. Hence
\[
 h\sqrt{1-q(s)}+\sum_r d_ru_r
 \le\sqrt{h^2+\sum_r u_r^2}<|F_{0,s}(z)|.
\]
Subtracting the rough envelope proves (OR2). At $s=0$ its right side is
$h\sqrt{w_*}=\delta$. Actual rough prime monomials are among the
independent multipliers allowed in (OR2), which proves the full-polydisc
conclusion.
\end{proof}

For each current source, all 53 inverse coefficients and all coefficients
of the product with the 1160-term smooth fibre are retained. Its integer
support has 3323 elements, including the constant. The complete defect
is less than $3/1000$. Thus no assertion that a truncated formal inverse
is itself zero-free is needed.

For the numerical recovery, let $v_r=Vw_r$ and partition all 157939
nonunit rough labels into 224 bins. In bin $k$ put
$m_k=\min r$ and $B_k=\sum_{r\ {\rm in\ bin}\ k}v_r$. Then
\[
 q(s)\le\widehat q(s)=\frac1V\sum_k B_km_k^{-2s}<1,
 \qquad \widehat R(s)=h\sqrt{1-\widehat q(s)}\le R(s).       \tag{OR3}
\]
Both $R$ and $\widehat R$ are nondecreasing. The equivalent formula
\[
 \widehat R(s)^2
 =\delta^2\left(1+\frac1{v_*}\sum_k B_k(1-m_k^{-2s})\right) \tag{OR4}
\]
makes positivity and $\widehat R(0)=\delta$ explicit, and provides an
independent numerical evaluation without subtracting two nearly equal
total weights. It follows by substituting $V=v_*+\sum_kB_k$ in (OR3).
All coefficient and bin uncertainties are enclosed outwards.

\subsection{Forward heat with a fixed finite support}

\begin{lemma}[Preservation of a closed-polydisc modulus]\label{lem:multi-heat}
Let $P$ be a finite polynomial nonzero on the closed unit polydisc. For
positive weights $\omega_j$ set
$D=\max_{\alpha\in\operatorname{supp}P}\omega\cdot\alpha$.
If $h\ge0$ and $4c\ge D$, then
\[
 \min_{\overline{\mathbb D}^{\,d}}
 \left|\sum_\alpha P_\alpha
 e^{h((\omega\cdot\alpha)^2/4-c\,\omega\cdot\alpha)}z^\alpha\right|
 \ge\min_{\overline{\mathbb D}^{\,d}}|P|.                   \tag{MA8}
\]
\end{lemma}

\begin{proof}
First consider $p_h=e^{h(\mathcal D^2/4-c\mathcal D)}p$ in one variable,
where $\mathcal D=z\partial_z$ and $p$ has degree $n$. Constant
polynomials are immediate. At a simple nonzero root $r_i$,
implicit differentiation and the factorisation of $p_h$ give
\[
 \frac{d}{dh}\log|r_i|
 =c-\frac14-\frac12\sum_{j\ne i}\Re\frac{r_i}{r_i-r_j}.
\]
If $|r_i|$ is minimal, $v=r_j/r_i$ has $|v|\ge1$ and $v\ne1$;
the identity
\[
 \Re\frac1{1-v}-\frac12=\frac{1-|v|^2}{2|1-v|^2}\le0
\]
therefore yields $d\log|r_i|/dh\ge c-n/4\ge0$.
On an interval of simple roots the minimum of their logarithmic radii is
locally Lipschitz and has nonnegative derivative almost everywhere, even
when the identity of the minimum changes. It is consequently nondecreasing.

If the initial roots are distinct, the discriminant is an analytic finite
exponential sum and is nonzero at $h=0$. It has finitely many zeros on
each compact real time interval. Root-multiset continuity joins the previous
conclusion across those collision times. Repeated initial roots are handled
by approximating them with distinct roots outside a common radius greater
than one. The constant and leading coefficients stay nonzero at finite
times. Taking the coefficient limit preserves nonvanishing in the closed
unit disk because the roots retain that common exterior radius.

Let $m=\min_{|z|\le1}|p(z)|>0$. Apply this stability assertion to
$p-\lambda$ for every $|\lambda|<m$. Constants are fixed by the
multiplier, so $p_h$ omits the whole disk $|\lambda|<m$ and
$|p_h|\ge m$ on the closed unit disk.

For positive integer weights $k_j$ and a fixed polydisc point $z$, apply
the univariate result to $p(\zeta)=P(z_j\zeta^{k_j})$. Its degree is
at most $\max k\cdot\alpha$, including any cancellations, and its
minimum modulus is at least the full-polydisc minimum of $P$.
Evaluation at $\zeta=1$ gives the desired bound at $z$.
For rational weights $\omega_j=k_j/Q$, use time $h/Q^2$ and drift $cQ$.
Finally choose positive rational vectors $\omega^{(m)}\to\omega$, let
$D_m=\max\omega^{(m)}\cdot\alpha$, and use
$c_m=c+\max(0,D_m-D)/4$. The corresponding finite coefficient sums
converge uniformly on the closed polydisc, proving (MA8), including
the case $4c=D$.
\end{proof}

\subsection{Closed moving-index gluing}\label{sec:multianchor}\label{sec:twoanchor}

For a finite source ladder let $C=J_1<J_2<\cdots<J_m<H$,
$J_{m+1}=H$, $Y_1=X$ and $Y_j=x_{J_j}(b)$ for $j\ge2$.
The four wall inequalities are
\[
 x_C(a)<X<x_{C+1}(b),\qquad x_H(a)<W<x_{H+1}(b).             \tag{MA9}
\]
They imply $N(t,X)=C$, $N(t,W)=H$ and $C\le N(t,x)\le H$
for $a\le t\le b$, $X\le x\le W$. They also put all later source
coordinates strictly between $X$ and $W$.

Cover $[J_j,J_{j+1}]$ by consecutive closed integer-endpoint leaves
$[L,U]$. On each leaf set
\[
 Z_L=\begin{cases}X,&j=1,\ L=C,\\x_{L}(b),&\text{otherwise},\end{cases}
 \qquad s_L=a\bigl(c_{\rm pol}(Z_L)-c_{\rm pol}(Y_j)\bigr),
\]
where in the second line $x_{L}(b)=4\pi(L^2-b/16)$; the symbol $Z_L$
on the left denotes the leaf's lower spatial coordinate. Define
\[
 f_L(v)=\exp\!\left\{\frac a4\log^2v-
                  \left(\frac12+a c_{\rm pol}(Z_L)\right)\log v\right\},
 \qquad I_{L,U}=\int_C^U f_L(v)\,dv.                       \tag{MA10}
\]

\begin{theorem}[Finite-exterior gluing]\label{thm:multi-gluing}
Assume every source satisfies (OR1) at its exact parameters, and every
leaf satisfies
\[
 \frac a2\log U-\frac12-a c_{\rm pol}(Z_L)<0,\qquad
 4c_{\rm pol}(Z_L)>\log U,\qquad
 \widehat R_j(s_L)-I_{L,U}>\mu>0.                           \tag{MA11}
\]
Then
\[
 \min_{\overline{\mathbb D}^{\,d_{N(t,x)}}}|A_{t,x,N(t,x)}|>\mu
 \quad (a\le t\le b,\ X\le x\le W).                       \tag{MA12}
\]
\end{theorem}

\begin{proof}
Fix the target $(t,x)$ and its actual right natural index $N=N(t,x)$.
Select a source interval containing $N$ and a leaf $[L,U]$ containing
$N$. Except on the first source's first leaf,
\[
 x\ge x_N(t)\ge x_{L}(b)=Z_L\ge Y_j.
\]
On that exceptional leaf $x\ge X=Z_L=Y_1$ is a domain hypothesis.
Also
\[
 c_{\rm pol}'(x)=\frac{x}{4(1+x^2)}+\frac{x}{(1+x^2)^2}>0.
\]
Thus $s=a(c_{\rm pol}(x)-c_{\rm pol}(Y_j))\ge s_L\ge0$, and
at the lower time the old coefficients obey the exact identity
\[
 a_n(a,x)=a_n(a,Y_j)n^{-s}\quad(n\le C).
\]
Ordinary recovery and its monotonicity give the lower bound
$\widehat R_j(s_L)$ for the complete old-support polynomial, throughout
the full prime polydisc. Extending this polynomial by unused new prime
variables does not alter that bound.

Every new coefficient belongs to the single set $C<n\le N$.
This includes smooth entrants, new rough fibres and coefficients between
$C$ and the selected source's own natural index. Spatial monotonicity
gives $a_n(a,x)\le f_L(n)$. Moreover
\[
 \frac{d}{dv}\log f_L(v)
 =\frac{a\log v/2-1/2-a c_{\rm pol}(Z_L)}{v}<0
 \quad(C\le v\le U)
\]
by the first inequality in (MA11). Therefore
\[
 \sum_{C<n\le N}a_n(a,x)
 \le\sum_{n=C+1}^U f_L(n)
 \le\sum_{n=C+1}^U\int_{n-1}^n f_L(v)\,dv
 =I_{L,U}.
\]
Since every entrant monomial has modulus at most one, reverse triangle
proves a strict lower bound $\widehat R_j(s_L)-I_{L,U}>\mu$ for
the full polynomial through the fixed target index $N$ at time $a$.

Keep $x$ and this complete support fixed while increasing the time from
$a$ to $t$. The multiplier on the $n$th coefficient is exactly
\[
 \exp\!\left\{(t-a)\left(\frac14\log^2n-c_{\rm pol}(x)\log n\right)\right\}.
\]
Apply Lemma~\ref{lem:multi-heat} with $\omega_p=\log p$. Its weighted
support degree is $\log N$, and
$4c_{\rm pol}(x)\ge4c_{\rm pol}(Z_L)>\log U\ge\log N$.
The minimum modulus is preserved. This proves (MA12).
The earlier support need not be natural at time $a$: it was deliberately
enlarged to the target support and all of its entrants were paid before
heat evolution.
\end{proof}

\subsubsection{Source handoffs and exact walls}

For a fixed time the natural cell of $N=n$ is
$[x_n(t),x_{n+1}(t))$. An index leaf $[L,U]$ includes every such
cell for each integer $L\le n\le U$. Consecutive closed leaves therefore
overlap on the whole common-index cell, and source intervals
$[J_{j-1},J_j]$ and $[J_j,J_{j+1}]$ overlap on the whole cell $N=J_j$.
No argument requires the two source formulae to agree there: each gives
a lower bound for the same complete target polynomial.

The fixed source coordinate $Y_j=x_{J_j}(b)$ lies to the left of the
moving lower cell wall when $t<b$. Explicitly,
\[
 \frac{Y_j}{4\pi}+\frac t{16}=J_j^2-\frac{b-t}{16},\qquad
 x_{J_j}(t)-Y_j=\frac\pi4(b-t).
\]
Here $0\le(b-t)/16\le43179/16000000<1$, so at $t<b$ the entire
thin interval $Y_j\le x<x_{J_j}(t)$ has $N=J_j-1$. It is covered
by the preceding source. At $t=b$, $N(b,Y_j)=J_j$ and either adjacent
source covers the point. At the moving wall itself the right index is
$J_j$, already covered. In particular
$N(a,Y_j)=J_j-1>C$ does not license an uncharged enlargement of the
source cutoff.

At any other index wall the same reasoning applies: the exact right
polynomial contains the entering term and has already paid it in the
sum from $C+1$ to the target index. This is a pointwise polynomial
statement. To infer a statement about the physical derivative at that
wall, the separate real-gauge theorem must hold for the frozen right
index on its local analytic disk, including the crossing remainder.
No derivative of the floor function is used in either step.

\subsubsection{Current certificate and the physical interface}

The 56 ordinary closed leaves give the following conservative rational
lower margins:
\[
\begin{array}{c|r|c|c}
j&\hbox{leaves}&\hbox{strict common floor}&\hbox{weakest leaf}\\ \hline
1&15&0.001272967812272&[690988,690995]\\
2&5&0.002111607661648&[693000,693875]\\
3&5&0.001597440849477&[730000,740000]\\
4&12&0.001368843821548&[791551,816691]\\
5&19&0.002893065121426&[1394929,1495493].
\end{array}                                                    \tag{MA13}
\]
Each terminating decimal is a strict rational lower bound and exceeds
$1/790$. No constrained recovery branch is used by this certificate.
For clarity, the actual whole-leaf acceptance test for every integer
index leaf $[L,U]$ is the interval inequality
\[
 \boxed{
   \underline{\widehat R_j(s_L)}
       - \overline{I_{L,U}}>\frac1{790},\qquad
   \overline{\frac{a\log U}{2}-\frac12-a c_{\rm pol}(Z_L)}<0,
   \quad \underline{4c_{\rm pol}(Z_L)-\log U}>0.}
                                                        \tag{RL}
\]
Here bars denote the outward upper and lower endpoints of directed
Arb enclosures, $s_L=a(c_{\rm pol}(Z_L)-c_{\rm pol}(Y_j))$ with
the stated exceptional first leaf, and $\widehat R_j$ is evaluated
from the cancellation-reduced formula \textup{(OR4)} with the certified
$P_{1,j},M_{1,j}$ and all 224 rough-weight bins. The entrant integral
is the complete integral from $C$ through $U$ below; no integer is
omitted. The 56 closed leaves, their stored enclosures, the source-bank
input hashes and an alternative 256-bit reconstruction are listed in
the reproducibility archive. The inequality is a leaf-wide statement,
not an evaluation at its centre.
For $B=1/2+a c_{\rm pol}(Z_L)-1$, the entrant integral is evaluated with
the exact primitive
\[
 \frac{\sqrt\pi}{\sqrt a}e^{-B^2/a}
 \operatorname{erfi}\!\left(\frac{a\log v/2-B}{\sqrt a}\right).
\]
Its derivative is $f_L(v)$; all values are enclosed outwards.

Combining the matching source Gram certificates, the full inverse defects
and (MA11) yields
\[
 |A_{t,x,N(t,x)}|>\frac1{790}
 \quad\left(\frac{156821}{1000000}\le t\le\frac15,
                 \ X\le x\le72382294738709\right)            \tag{FE}
\]
throughout the closed prime polydisc. This gives both stability and the
torus modulus premise of Corollary~\ref{cor:admission}. The exact
coefficient identity at the beginning of this section ensures that the
same natural-index polynomial enters that corollary. Since
$a>98/625$ and $W<120000000000000$, its parameter domain contains
the whole finite exterior, including both endpoints and every index wall.
The corollary consequently excludes physical multiple real zeros there.
The stationary separator, interior argument, infinite exterior and global
strictness argument remain separate components of the global conclusion.

\section{Validated polynomial inequalities: arithmetic and enclosure criteria}

For each exact source, all coefficients of the real trigonometric
polynomial $T_j$ in (G) are reconstructed by directed Arb arithmetic.
The finite support is generated from integer prime-valuation vectors,
with the full rough-fibre pair sums and the full reciprocal-product
convolution. A Fourier coefficient outside this structural support is
algebraically zero. There are 132901 represented frequencies per source,
including opposite frequencies and the constant.

Write $T=\sum_k c_ke^{ik\cdot\theta}$ and split it as a retained
finite real cosine polynomial $T_{\rm keep}$ plus a remainder. The
outward sum of the absolute values of all discarded coefficients is
$\eta$, so
\[
 |T(\theta)-T_{\rm keep}(\theta)|\le\eta
 \quad\hbox{at every phase.}
\]
Both coefficient enclosure uncertainty and trigonometric quantisation
are included in the derivative-jet enclosures. In particular a discarded
small coefficient is charged uniformly, rather than set to zero without
an error allowance.

For a closed box with centre $\theta_0$ and half-widths
$h=(h_1,\ldots,h_4)$, let $J_\alpha$ enclose
$\partial^\alpha T_{\rm keep}(\theta_0)$ for $|\alpha|\le3$,
and let
\[
 M_\alpha\ge\sum_{k\ \rm kept}|c_k|\prod_{r=1}^4|k_r|^{\alpha_r}
 \qquad(|\alpha|=4).
\]
The cosine representation collects the opposite modes with their correct
factor of two. The displayed formula is equivalently stated in the full
Fourier representation. Taylor's formula along the segment from
$\theta_0$ to $\theta$ gives the rigorous lower bound
\[
 T(\theta)\ge\inf J_0-\eta
 -\sum_{1\le|\alpha|\le3}\frac{h^\alpha}{\alpha!}
          \sup|J_\alpha|
 -\sum_{|\alpha|=4}\frac{h^\alpha}{\alpha!}M_\alpha.       \tag{FB1}
\]
Indeed the fourth derivative of that one-variable restriction expands
multinomially, and $|e^{ik\cdot\theta}|=1$ bounds its every derivative
term by the corresponding absolute Fourier moment. The one-variable
integral remainder supplies the factor $1/4!$, cancelling the
multinomial numerator. This proves (FB1) on the closed box, including
all its faces. There are 35 multiindices through degree three and 35
of degree four in four variables.

The programs quantise the retained coefficients and trigonometric
values to denominator $2^{44}$, compute the jets by exact integer
accumulation, and bound every quantisation error before (FB1) is used.
Signed 64-bit input and phase ranges and signed 128-bit accumulator
ranges are verified before the C++ calculation. A Python arbitrary-
precision implementation checks selected complete rows, including
every maximum-depth row and each critical row. The range argument and
the full native loop apply to every row; the selected comparisons are
additional implementation controls.

Each final box comes from a finite exact dyadic cover. The primary
coverage check establishes total volume and disjoint interiors. A
separate check starts from an exact dense census of the 6016-box
initial cover on $64^4$ cells, then replays every bisection and demands
equality with the final box list. Closed shared faces belong to both
neighbours. Thus checking (FB1) on all final boxes proves the whole
torus inequality, with no appeal to an unexamined numerical sampling
density. Every final lower bound exceeds the exact threshold $10^{-6}$.

An additional export audit inspects every dyadic coefficient interval
$[l_k/2^{80},u_k/2^{80}]$ and proves $u_k-l_k\le2$.
Consequently each exported midpoint has error at most $2^{-80}$.
This width conclusion is checked from the actual intervals; it does
not follow merely from the choice of an export grid.
The original $0.16$ comparison separately includes the absolute mass
of floating coefficients outside the algebraic support, so its full
Fourier $\ell^1$ discrepancy also controls the uniform polynomial error.

\section{A stationary nonvanishing separator on the initial time band}\label{sec:separator}

The stationary-separator architecture draws on the pinned computational
record \cite{MosaicAssembly}. This section supplies the physical
transfer for the independently reconstructed moments used here.

Set
\[
 \begin{gathered}
 X=6000000185827,\quad N=690988,\quad n_0=N/2,\quad T=1/5,\\
 Q=X/(4\pi),\quad L=\log Q,\quad \ell=\log n_0.
 \end{gathered}
\]
The following is a direct application of Polymath's explicit Theorem 1.3,
equations (20)--(24), with the error bounds retained uniformly.  Its time
parameter is within the published domain $0<t\le1/2$; no unpublished
starting bound is involved.

For $0\le y\le1$ write $s=(1-y+iX)/2$, $h=1-s$ and
\[
 \alpha(w)=\frac1{2w}+\frac1{w-1}+\frac12\Log\frac{w}{2\pi},
 \qquad \gamma=\frac{M_t(s)}{M_t(h)},\qquad B_t=M_t(h).
\]
Here the elementary prefactor is explicitly
\[
 M_0(w)=\frac{\sqrt{2\pi}}{16}w(w-1)
 \exp\left[-\frac w2\log\pi+
 \left(\frac w2-\frac12\right)\Log\frac w2-\frac w2\right],
 \qquad M_t(w)=M_0(w)e^{t\alpha(w)^2/4}.
\]
All logarithms are their principal branches.  The arguments stay off the cuts,
and $B_t\ne0$: $\Im h=-X/2\ne0$ makes $h(h-1)$ nonzero, and every remaining
factor in $M_t(h)$ is a nonzero constant or exponential.  Since
$s(s-1)=h(h-1)$ exactly, cancellation gives the formula used in the code,
\[
 \gamma=\exp\left[-\frac{s-h}{2}(\log\pi+1)
 +\left(\frac s2-\frac12\right)\Log\frac s2
 -\left(\frac h2-\frac12\right)\Log\frac h2
 +\frac t4\bigl(\alpha(s)^2-\alpha(h)^2\bigr)\right].
\]

Define
\[
 D_t(w)=\sum_{n=1}^{N}e^{t(\log n)^2/4}n^{-w-t\alpha(w)/2},
 \qquad f_t=D_t(h)+\gamma D_t(s).
\]
To identify this with the published finite model, let
$s_*=h+t\alpha(h)/2$ and
$\kappa=t(\alpha(s)-\alpha(\overline h))/2$.  Then
$\overline{s_*}+\kappa-y=s+t\alpha(s)/2$ exactly.  Thus the second
published chirality is precisely $D_t(s)$; no conjugation or heat-phase
correction is dropped.

The exact natural-index gates
\[
 N^2<Q,\qquad Q+T/16<(N+1)^2
\]
hold with respective lower margins exceeding $427919$ and $954057$.
Thus the same $N$ is the natural index on the entire closed time interval.

\subsection{Exact polynomial reconstruction and its truncation}

Put $\theta=(2X+1)/4$ and $L_n=\log(n/n_0)$.  The fresh Acb source bank
contains enclosures for the exact sums
\[
 S_k=\sum_{n=1}^{N}n^{-1/2}e^{i\theta\log n}L_n^k,
 \qquad 0\le k\le183.
\]
With $C_{e,j}=S_{e+2j}/(e!j!4^j)$ define
\[
 P_t(z)=\sum_{e=0}^{61}\sum_{j=0}^{61}C_{e,j}z^e t^j.
\]
Let
\[
 b_{\rm sep}=-y/2-i/4-t\alpha(h)/2,\qquad
 a_{\rm sep}=y/2-i/4-t\alpha(\overline s)/2,
\]
\[
 \widetilde B=e^{(b_{\rm sep}+t\ell/4)\ell}P_t(b_{\rm sep}+t\ell/2),\qquad
 \widetilde A=e^{(a_{\rm sep}+t\ell/4)\ell}P_t(a_{\rm sep}+t\ell/2),\qquad
 \widetilde f=\widetilde B+\gamma\overline{\widetilde A}.
\]
These identities are obtained by writing $\log n=\ell+L_n$ in each chirality.
The quarter shift is exactly $X/2-\theta=-1/4$.  Replacing each truncated
series by its full exponential series recovers $D_t(h)$ and
$D_t(\overline s)$ exactly; conjugation supplies $D_t(s)$.

Whole-domain Arb evaluation proves
\[
 |a_{\rm sep}+t\ell/2|<2/3,\quad |b_{\rm sep}+t\ell/2|<2/3,\quad
 \Re\alpha(h)/2>\ell/4,\quad\Re\alpha(\overline s)/2>\ell/4.
\]
Thus the two external factors have moduli at most $\sqrt{n_0}$ and $1$,
respectively.  Published equation (20) gives
$|\gamma|\le e^{y/50}Q^{-y/2}\le1$, since $L/2>1/50$.
Also $|L_n|\le\ell$ for every integer $1\le n\le2n_0$.

For $T_{61}(u)=\sum_{j=0}^{61}u^j/j!$, the elementary exponential
remainder and $|T_{61}(u)|\le e^{|u|}$ give
\[
 |e^{u+v}-T_{61}(u)T_{61}(v)|
 \le e^{|u|+|v|}\frac{|u|^{62}+|v|^{62}}{62!}.
\]
Using $U=2\ell/3$, $V=T\ell^2/4$ and $n^{-1/2}\le1$ therefore yields
\[
 |f_t-\widetilde f_t|\le E_T
 :=N(1+\sqrt{n_0})e^{U+V}\frac{U^{62}+V^{62}}{62!}
 <9.811176\,10^{-13}.
\]
All moment and coefficient enclosure radii are already retained by the
whole-box Acb evaluation; this is a separate Taylor truncation charge.

\subsection{Published approximation charges, uniformly in $t,y$}

Define the exact positive parameters
\[
 \delta=\frac{T(1+8/X^2)}{2X^2},\qquad
 k=\frac{T}{2(X-6)},\qquad
 u=\frac{T^2L^2/16+313/500}{X-333/50}.
\]
Published equation (21) bounds the real saddle exponent from below by
$(1+y)/2+tL/4-\delta$: the positive part of
$1-3y+4y(1+y)/X^2$ is at most $1+8/X^2$ on $[0,1]$.
The coefficient factor
$e^{t((\log n)^2-L\log n)/4}$ is at most one.  Equations (20) and (22),
together with $n\le N\le\sqrt Q$, bound the second-to-first chirality factor
by $e^{1/50+k\log N}$.  Since
$0\le\log(Q/n^2)\le L$ and $\sum_{n\le N}n^{-1/2}\le2\sqrt N$,
equation (23) gives
\[
 e_A+e_B\le E_{AB}:=
 2\sqrt N\bigl(1+e^{1/50+k\log N}\bigr)
 e^{\delta\log N}\operatorname{expm1}(u)
 <1.362431\,10^{-9}.
\]
Dropping only nonpositive exponents in equation (24), and using
$3^y+3^{-y}\le4$, gives
\[
 e_{C,0}\le E_C:=Q^{-1/4}
 \exp\left\{\frac{(31/25)4}{N-1/8}
 +\frac{3\sqrt{L^2+(\pi/2)^2}+261/25}{X-12}\right\}
 <0.001203005928212.
\]
All denominators are positive.  Every time-dependent upper envelope uses
$T=1/5$, not the earlier $3/16$.  Fresh 256-bit Arb evaluation proves
\[
 E_T+E_{AB}+E_C
 <0.001203007291623<1/800.
\]
The approximation theorem is stated for positive $t$.  The fixed-$N$ model,
$H_t$ and the nonzero prefactor are continuous as $t\downarrow0$.  The common
uniform error bound therefore extends to $t=0$ as a non-strict error
inequality, which suffices for the strict final reserve.

\subsection{New closed-box witness and conclusion}

The fresh certificate covers every point of
$[0,1/5]\times[0,1]$ with $17385$ closed rectangles, obtained from $10250$
exact rational seeds by $7135$ exact bisections.  Every whole-box ball
lower bound is strictly above $1/100$, with common valid floor
\[
 |\widetilde f_t(X+iy)|>0.0100458444330390.
\]
The largest depth is two.  A separate rational geometry check verifies
parent containment, exact area, absence of positive-area overlaps and all
seed cells; finite closedness includes seams and corners.
For each closed rectangle $R=I_t\times I_y$ (both intervals have
rational endpoints), the certification program forms an Acb enclosure
\[
 \mathcal F(R)=\widetilde B(I_t,I_y)
       +\gamma(I_t,I_y)\,\overline{\widetilde A(I_t,I_y)}
       \quad\text{with}\quad
 \{\widetilde f_t(X+iy):(t,y)\in R\}\subseteq\mathcal F(R),
\]
using the 184 complex moment enclosures, and accepts $R$ only if
\[
 \boxed{\operatorname{dist}(0,\mathcal F(R))
       \ge \operatorname{abs\_lower}(\mathcal F(R))>\frac1{100}.}
                                                        \tag{SR}
\]
The program uses the rectangular interval substitutions $I_t,I_y$
\emph{throughout} its evaluation and calls the direct Acb
\texttt{abs\_lower()} method; it does not sample only the vertices or
centres. Its final CSV records the rational rectangle and outward
lower bound, while the separate geometry checker verifies the exact
closed cover. This makes the universal meaning of every reported
$17385$-rectangle assertion explicit.  The previously
uncovered strip $3/16\le t\le1/5$ has $732$ final rectangles.

Subtracting the error charges proves
\[
 \boxed{\left|\frac{H_t(X+iy)}{B_t(X+iy)}\right|
 >0.0088428371414166>\frac7{800}}
 \qquad(0\le t\le1/5,\ 0\le y\le1).
\]
Because $B_t$ is nonzero, evenness and conjugation supply both stationary
vertical boundaries, for every $-1\le y\le1$.  Together with finite RH
initialisation and Theorem~\ref{thm:compact-interior}, this proves that all
zeros with $|\Re z|\le X$ are real and simple for every $0<t\le1/5$.

Platt--Trudgian's Corollary 2 supplies the published global anchor
$\Lambda\le1/5$. The exterior and practical tail proved in
Sections~\ref{sec:source}--\ref{sec:tail} cover the same enlarged upper
time endpoint. Section~\ref{sec:descent} combines them with this separator
and the following interior theorem to prove the strict global conclusion.

\subsection{Complex-radius correction and the common normalisation}

The supplementary floating-moment audit in the original Planat package
assigns a radius $\rho_k$ to each of the real and imaginary coordinates.
Such a rectangular enclosure gives the scalar estimate
\[
 |\Delta S_k|\le\sqrt{(\Delta\Re S_k)^2+(\Delta\Im S_k)^2}
              \le\sqrt2\,\rho_k.
\]
Accordingly its displayed scalar error allowance $4.211\cdot10^{-7}$
must be replaced by the conservative bound $5.96\cdot10^{-7}$ unless a
separate disk enclosure of the smaller radius is proved. This correction
concerns that supplementary scalar ledger. The original direct complex
interval box calculations already propagate both coordinate intervals.
The newly rebuilt Acb calculation likewise propagates both coordinates
through every operation. An additional $\sqrt2$ debit on those complete
complex enclosures would count the same uncertainty twice.

On the common real line the balanced gauge satisfies
\[
 M_t(s_+)=GP,\qquad B_t=M_t(s_-)=G/P,\qquad
 \frac{H_t}{B_t}=Pu,\qquad
 \partial_x\left(\frac{H_t}{B_t}\right)
       =P\bigl(u_x+\ii\Theta_xu\bigr).
\]
Here $G>0$, $\Theta$ is real and $|P|=1$. Only the first normalised
modulus is used for the stationary boundary. The finite and infinite
exterior use the real pair $(u,u_x)$, with the carrier derivative retained.
No equality of differently normalised off-axis moduli is required.

\section{Propagation of real zeros in the compact interior}
\label{sec:interior}

The argument below permits multiple zeros at the initial time. Increasing
the parameter means evolution by $\partial_t F=-\partial_z^2F$ throughout
this section; this sign convention is essential.

\begin{lemma}[Local splitting at a real zero]\label{lem:local-heat-splitting}
Let $F$ be jointly holomorphic near $(t_0,x_0)\in\mathbb R^2$, satisfy
$\partial_tF=-\partial_z^2F$, and be real-type for real $t$. Suppose
$F(t_0,\cdot)$ has a zero of finite multiplicity $m\ge1$ at $x_0$.
In a sufficiently small fixed disk about $x_0$, every zero of
$F(t_0+h,\cdot)$ is real and simple for all sufficiently small $h>0$.
There are exactly $m$ such zeros, counted with multiplicity. If $m\ge2$,
then $F(t_0-h,\cdot)$ has a nonreal zero approaching $x_0$ for every
sufficiently small $h>0$.
\end{lemma}

\begin{proof}
Write $c=\partial_z^mF(t_0,x_0)/m!\ne0$. The differential equation gives
\[
 \partial_t^j\partial_z^k F(t_0,x_0)
       =(-1)^j\partial_z^{2j+k}F(t_0,x_0).
\]
Consequently the joint Taylor coefficient of
$(t-t_0)^j(z-x_0)^k$ vanishes whenever $2j+k<m$.
For a real parameter $s>0$, substitute
$(t,z)=(t_0\pm s^2,x_0+sw)$. The resulting holomorphic functions of
$(s,w)$ are divisible by $s^m$. Their quotients extend holomorphically
to $s=0$, locally uniformly for $w$ in the complex plane, and satisfy
\begin{align}
 s^{-m}F(t_0+s^2,x_0+sw)&=cP_m(w)+O(s),\label{eq:hermite-forward}\\
 s^{-m}F(t_0-s^2,x_0+sw)&=cQ_m(w)+O(s),\label{eq:hermite-backward}
\end{align}
uniformly on every compact $w$-set. Here
\[
 P_m(w)=\sum_{j=0}^{\lfloor m/2\rfloor}
 \frac{(-1)^j m!}{j!(m-2j)!}w^{m-2j},\qquad
 Q_m(w)=\sum_{j=0}^{\lfloor m/2\rfloor}
 \frac{m!}{j!(m-2j)!}w^{m-2j}=i^{-m}P_m(iw).
\]
This substitution also proves the uniform error assertion: on a compact
$w$-set, the quotient is holomorphic for $|s|$ sufficiently small, and its
Taylor series at $s=0$ has the displayed constant term.

The polynomial $P_m$ is monic and has $m$ distinct real roots. A short
self-contained proof follows from its generating function
\[
 \sum_{m\ge0}P_m(w)\frac{u^m}{m!}=e^{uw-u^2}
       =e^{w^2/4}e^{-(w-2u)^2/4}.
\]
Thus
\[
 P_m(w)=(-2)^m e^{w^2/4}
            \frac{d^m}{dw^m}e^{-w^2/4}.
\]
For every real polynomial $q$ of degree less than $m$, integration by parts
$m$ times, with all boundary terms zero, gives
\[
 \int_{\mathbb R}P_m(w)q(w)e^{-w^2/4}\,dw=0.
\]
If $P_m$ had fewer than $m$ distinct real roots of odd multiplicity, take
$q$ to be the product of the corresponding linear factors. The product
$P_mq$ then has constant sign on $\mathbb R$ and is not identically zero,
contradicting this identity. Hence all $m$ roots are real and simple.
Parity makes the roots symmetric about zero. For $m\ge2$, at least one
pair is nonzero, so $Q_m$ has a nonreal pair of simple roots; its roots
are purely imaginary, with a possible root at zero when $m$ is odd.

Choose disjoint disks centred at the real roots of $P_m$. Each disk is
invariant under complex conjugation. Rouch\'e's theorem and
\eqref{eq:hermite-forward} put exactly one zero, counted with multiplicity,
in each of their images under $w\mapsto x_0+sw$, for every sufficiently
small $s>0$. Real-type symmetry forces this unique zero to be real; its
count one implies simplicity. Separately, choose a fixed disk about
$x_0$ containing no other zero of $F(t_0,\cdot)$ and no zero on its
boundary. Rouch\'e's theorem on this fixed boundary gives exactly $m$
zeros at nearby times. For small $s$, all the scaled disks lie inside it,
so their $m$ zeros exhaust the cluster.

For $m\ge2$, apply \eqref{eq:hermite-backward} on a disk about a nonreal
root of $Q_m$ whose closure misses the real axis. Its image contains a
nonreal zero of $F(t_0-s^2,\cdot)$ for every sufficiently small $s>0$.
Its distance from $x_0$ tends to zero. This proves the lemma with $h=s^2$.
\end{proof}

\begin{theorem}[Compact interior propagation]\label{thm:compact-interior}
Let $T,X,Y>0$ and
\[
 \Omega=\{z:|\Re z|<X,\ |\Im z|<Y\}.
\]
Suppose $F$ is jointly holomorphic on a complex neighbourhood of
$[0,T]\times\overline\Omega$, satisfies
$\partial_tF=-\partial_z^2F$, and is real-type for real $t$.
Assume
\[
 F(t,z)\ne0\quad(0\le t\le T,\ z\in\partial\Omega),
\]
and that every zero of $F(0,\cdot)$ in $\Omega$ is real.
Then every zero of $F(t,\cdot)$ in $\Omega$ is real and simple for
$0<t\le T$. No hypothesis of initial simplicity is required.
\end{theorem}

\begin{proof}
The boundary condition excludes the identically zero function at every
time. It also gives a positive minimum of $|F|$ on the compact set
$[0,T]\times\partial\Omega$. By the argument principle the number $M$
of zeros in $\Omega$, counted with multiplicity, is finite and independent
of $t$.

At any time, take small disjoint disks about the finitely many interior
zeros, with closures in $\Omega$ and no zeros on their boundaries. The
closed complement in $\overline\Omega$ is zero-free. Uniform continuity
on this complement and Rouch\'e's theorem on the disk boundaries show
that, at all sufficiently nearby times, each disk contains the same
number of zeros and there are no other zeros in $\overline\Omega$.
This supplies the required continuity of clusters without labelling
individual zero trajectories.

If $M=0$ there is nothing to prove. Otherwise apply
Lemma~\ref{lem:local-heat-splitting} to the finitely many initial zeros.
Together with the preceding compact-complement argument, this shows
that all interior zeros are real and simple on some interval $0<t\le d$.

Let $\tau$ be the supremum of the numbers $v\in[d,T]$ for which this
property holds at every $0<t\le v$. The assertion is already proved if
$d\ge T$. Otherwise $0<\tau\le T$. Every zero at time $\tau$ is real:
an off-real zero would persist in a disk missing the real axis at times
just below $\tau$, which is impossible. No such zero can be multiple,
because the backward part of Lemma~\ref{lem:local-heat-splitting} would
produce a nonreal zero at every sufficiently close earlier time. These
earlier zeros remain in $\Omega$ since their limiting zero is interior.

All zeros at $\tau$ are therefore real and simple. The cluster-continuity
argument now uses conjugation-invariant disks containing one zero each,
so it preserves reality and simplicity for a neighbourhood of $\tau$.
If $\tau<T$ this contradicts the definition of $\tau$; if $\tau=T$ it
also proves the result at the final endpoint. Thus the conclusion holds
on the whole claimed interval.
\end{proof}

\subsection{Application to the Riemann $\xi$ flow}
Use
\[
 H_t(z)=\int_0^\infty e^{tu^2}\Phi(u)\cos(zu)\,du,
 \qquad H_0(z)=\frac18\xi\left(\frac{1+iz}{2}\right).
\]
The super-exponential decay of $\Phi$ makes the integral and all its
parameter derivatives locally uniformly convergent in $(t,z)\in\mathbb C^2$.
Thus $H$ is jointly entire, is even and real-type when $t$ is real, and
satisfies $\partial_tH_t=-\partial_z^2H_t$.

The initial critical strip for $\xi$, together with de Bruijn's strip
theorem in the form \cite[Theorem~3.2]{Polymath2019}, gives
\[
 H_t(z)=0,\ t\ge0
 \quad\Longrightarrow\quad
 |\Im z|\le\sqrt{\max(1-2t,0)}\le1.
\]
Only the final coarse inequality is needed below. At $t=0$ it follows
directly from the critical strip; for $t>0$ it follows by applying the
cited theorem at initial time zero and using conjugation symmetry.

Set
\[
 b=\frac15,\qquad X=6000000185827,\qquad Y=2.
\]
The separator hypothesis required here is exactly
\begin{equation}\label{eq:closed-physical-separator}
 H_t(X+iy)\ne0
 \qquad(0\le t\le b,\ 0\le y\le1).
\end{equation}
Every endpoint in this statement is included. A certified nonzero value
of $H_t/B_t$ implies \eqref{eq:closed-physical-separator} provided that
the actual normalising factor $B_t$ is nonzero on this domain.

Conjugation reflects \eqref{eq:closed-physical-separator} across the real
axis, and evenness supplies the side $\Re z=-X$. The strip bound excludes
zeros on the remaining portions of both vertical sides, where
$1<|\Im z|\le2$, and on the horizontal sides $\Im z=\pm2$.
This proves the full closed-rectangle boundary condition, including all
corners, for every $t\in[0,b]$.

The affine zero-coordinate map is
\[
 \rho=\beta+i\gamma\quad\longmapsto\quad
 z=2\gamma+i(1-2\beta).
\]
It preserves multiplicities. Platt--Trudgian's exact finite-RH result
\cite[Theorem~1]{PlattTrudgian2021} reaches height
$T_{\rm PT}=3000175332800$, and
\[
 T_{\rm PT}-\frac X2=\frac{350479773}{2}>0.
\]
Hence every initial zero with $0<|\Re z|\le X$ is real, using the
conjugation symmetry of zeta for negative ordinates. There is no initial
zero with $\Re z=0$: for real $y$,
$H_0(iy)=\int_0^\infty\Phi(u)\cosh(yu)\,du>0$.
Indeed every term of $\Phi(u)$ for $u\ge0$ equals
$\pi n^2e^{5u}(2\pi n^2e^{4u}-3)e^{-\pi n^2e^{4u}}>0$.
The exact finite-height endpoint is needed; the rounded height in the
title of the cited paper would be slightly smaller than $X/2$.

\begin{corollary}[Physical finite interior]\label{cor:physical-interior}
If \eqref{eq:closed-physical-separator} holds, every zero of $H_t$ with
$|\Re z|\le X$ is real and simple for $0<t\le1/5$. There are no zeros
at either vertical endpoint $\Re z=\pm X$.
\end{corollary}
\begin{proof}
Apply Theorem~\ref{thm:compact-interior} with $F=H$, $T=1/5$ and $Y=2$.
The strip bound places every zero with $|\Re z|<X$ inside that rectangle.
The closed vertical endpoints are zero-free by the separator and strip
arguments just given.
\end{proof}

This corollary consumes the physical separator, not merely the polynomial
moment certificate. Its application therefore requires the complete
separator error budget and the $t=0$ limit of the chosen normalisation.
It consumes neither the $0.1732$ certificate nor the global starting
bound. The latter enters only the global assembly.

\section{Uniform exclusion of multiple zeros in the infinite exterior}
\label{sec:tail}

Throughout this section the constants are exactly
\[
 H_*=2400000,\qquad a=\frac{156821}{1000000},\qquad b=\frac15,
 \qquad W=72382294738709=\lceil4\pi H_*^2\rceil,
\]
\[
 \begin{gathered}
 \mathcal P=\{2,3,5,7,11,13\},\quad D=30030,\quad K=2391233,\\
 K_0=K-D=2361203,\quad\epsilon_0=H_*^{-2},\quad U_0=\log H_*.
 \end{gathered}
\]
This section proves real-collision exclusion for every $a\le t\le b$ and $x\ge W$.
The balanced-gauge theorem is used with its proved actual bounds
\[
 |u-V|<v=\frac{101}{10^8},\qquad
 |u_x-V_x|<d=\frac{71}{10^7}.
                                                               \tag{ET1}
\]
Here $u=H_t/G$ with $G(x)>0$, and $V$ is the two-chiral real model.
Those bounds hold at all $x\ge X=6000000185827$ on the larger time band
$[98/625,1/5]$, so they apply on the entire present tail.

\subsection{The natural wall and the exact saddle}

At the point being considered set
\[
 N=\left\lfloor\sqrt{x/(4\pi)+t/16}\right\rfloor,\qquad U=\log N.
\]
Use the right natural index at an exact wall, and hold that index fixed
for local spatial differentiation. Define
\[
 s=\frac{1+ix}{2},\qquad
 \alpha(s)=\frac1{2s}+\frac1{s-1}+\frac12\Log\frac{s}{2\pi},
 \qquad w=s+\frac t2\alpha(s)=\sigma+i\tau.
\]
The real part is exactly
\[
 \sigma=\frac12+tc_{\rm pol}(x),\qquad
 c_{\rm pol}(x)=\frac14\log\frac{\sqrt{1+x^2}}{4\pi}
                         -\frac1{2(1+x^2)}.
\]
Put $\bar\sigma(U,t)=1/2+t(U/2-\epsilon_0)$.

For the present domain $x\ge W$, the wall definition gives
$N\ge H_*$; thus this is the only case. The natural wall implies
$x\ge4\pi(N^2-t/16)\ge4\pi(N^2-b/16)>N^2$.
For $z=b/(16N^2)$, the inequality
$\log(1-z)\ge-z/(1-z)$ gives
\[
 \log\frac{x}{4\pi}\ge2U-\frac{b}{16N^2-b}
                  >2U-\frac1{N^2}\ge2U-\epsilon_0.             \tag{ET2}
\]
The rational endpoint gate $b/(16H_*^2-b)<H_*^{-2}$ verifies the strict
step, and the inequality persists as $N$ increases. Also
\[
 c_{\rm pol}(x)\ge\frac14\log\frac{x}{4\pi}-\frac1{2x^2}
       >\frac U2-\frac{\epsilon_0}{4}-\frac1{2H_*^4}
       >\frac U2-\epsilon_0,
\]
where $1/(4H_*^2)+1/(2H_*^4)<1/H_*^2$ is another exact rational gate.
Consequently $\sigma\ge\bar\sigma(U,t)$ throughout the tail.

For $b_n(t)=e^{t\log^2n/4}$, define
\[
 A_N(w)=\sum_{n\le N}b_n(t)n^{-w},\qquad
 E(w)=\prod_{p\in\mathcal P}(1-b_p(t)p^{-w}).
\]
For square-free $d\mid D$ put
\[
 q_d=\frac14\sum_{p\mid d}\log^2p,\qquad
 \lambda_d(t)=\mu(d)e^{tq_d}.
\]
Every represented product in $EA_N$ through $DN$ is retained. For
$m\le K<N$, its coefficient is the complete convolution
\[
 \beta_m(t)=\sum_{d\mid\gcd(m,D)}\lambda_d(t)b_{m/d}(t).
                                                               \tag{ET3}
\]
The constant coefficient is one.  Products above $K$, including those
above $N$, will be bounded by their complete individual pairs.

\subsection{Positivity and uniform comparison in heat time}

Let the selected primes dividing $m$ be vertices, with labels
$\ell_p=\log p$, and put
$R_m=\log(m/\prod_{p\mid\gcd(m,D)}p)\ge0$.
At time $t$, place independent Poisson arrivals on every pair of
distinct vertices at rate $\ell_p\ell_q/2$, and on the loop at vertex
$p$ at rate $R_m\ell_p/2$. An arrival covers its incident vertices.
For a subset $S$ of vertices, the rate of all edges incident to $S$ is
\[
 \frac12\log m\sum_{p\in S}\ell_p
 -\frac14\left(\sum_{p\in S}\ell_p\right)^2
 -\frac14\sum_{p\in S}\ell_p^2.
\]
Inclusion-exclusion for the event that no vertex remains uncovered
therefore gives exactly
\[
 \beta_m(t)=b_m(t)C_m(t)\ge0,                                  \tag{ET4}
\]
where $C_m(t)$ is the probability that all vertices are covered.
The term $R_m$ retains repeated prime powers and every unselected
prime factor.  If there are no selected vertices, $C_m=1$. If $m$ is
itself one selected prime, there is one vertex and no positive-rate
edge, so $\beta_m=0$ identically.

If there are $k$ vertices and a realised graph covers them, it contains
a covering subgraph with at most $k$ arrived edges: choose one incident
arrival for each vertex and discard duplicates.  Thin arrivals
independently from time $t_2$ to $t_1<t_2$. The chance of retaining those
chosen arrivals is at least $(t_1/t_2)^k$. Hence
\[
 C_m(t_1)\ge(t_1/t_2)^kC_m(t_2),\qquad
 \frac{C_m'(t)}{C_m(t)}\le\frac{k}{t}
\]
where $C_m(t)>0$.  For $m=6$, the single pair edge covers both vertices,
so the sharper $1/t$ applies.

For $f_U(v)=v(U/2-\epsilon_0-v/4)$ the logarithmic heat derivative of
the positive head term is at most $k/t-f_U(\log m)$.
The function $f_{U_0}$ is concave.  Fresh directed endpoint checks give
\[
 f_{U_0}(\log m_k)>k/a,\qquad f_{U_0}(\log K)>k/a
 \quad(k=1,\ldots,6),
\]
\[
 (m_1,m_2,m_3,m_4,m_5,m_6)=(4,10,30,210,2310,30030),\qquad
 f_{U_0}(\log6)>1/a.                                          \tag{ET5}
\]
These are the minimum possible nonzero indices for each $k$, with
$m=6$ treated separately.  Concavity proves the inequalities throughout
the respective complete ranges.  They persist for $U\ge U_0$ and
$t\ge a$, because $f_U$ increases with $U$ and $k/t\le k/a$.
In the no-vertex case $f_U(\log m)>0$ for $m\ge2$, since
$\log m\le\log K<U_0$ and $U_0/4>\epsilon_0$.
Thus every complete nonconstant positive head term
$\beta_m(t)m^{-\bar\sigma(U,t)}$ decreases with heat time.

For a remaining product pair $(d,n)$, $d\mid D$ and $n\le N$, its
amplitude at the majorising real part is
\[
 e^{tq_d-\bar\sigma(U,t)\log d}
 e^{t\log^2n/4-\bar\sigma(U,t)\log n}.
\]
It also decreases with $t$. Indeed, $q_d\le(\log d)^2/4$ and the
fresh signs
\[
 \log D/4<U_0/2-\epsilon_0,\qquad U_0/4>\epsilon_0              \tag{ET6}
\]
make both contributions to its logarithmic heat derivative nonpositive.
The logarithmic derivative of its $n$-density is at most
$-1/2+b\epsilon_0<0$. These comparisons include all products $dn>N$.

\subsection{A complete finite head and a continuous complement barrier}

Set
\[
 \mathcal H(U)=\sum_{m=2}^K\beta_m(a)m^{-\bar\sigma(U,a)},
 \qquad \kappa=\frac a2\log K_0,
\]
\[
 C_0=\sum_{2\le m\le K_0}
 (\kappa-\tfrac a2\log m)\beta_m(a)m^{-\bar\sigma(U_0,a)}\ge0.
\]
Differentiating each positive summand gives
$\mathcal H'\le-\kappa\mathcal H+C_0$: for $m>K_0$ the extra
coefficient is negative, and for $m\le K_0$ its positive contribution
decreases as $U$ increases.

For the remaining pairs, $n\ge\lfloor K/d\rfloor+1$.  The decreasing
density and (ET6) give the uniform upper bound
\[
 \mathcal R(U)=\sum_{d\mid D}e^{aq_d-\bar\sigma(U,a)\log d}
 \int_{\log\lfloor K/d\rfloor}^{U}
 e^{av^2/4-(\bar\sigma(U,a)-1)v}\,dv.                          \tag{ET7}
\]
The lower limit is the logarithm of the floor; replacing it by
$\log(K/d)$ or its ceiling would not give the stated upper bound.
Every integral is positive since $K<H_*\le N$.
For its internal derivative,
\[
 \log d+v\ge\log(d\lfloor K/d\rfloor)>\log(K-D)=\log K_0,
\]
so the derivative is at most $-\kappa$ times the integrand.
Its top-boundary forcing is
\[
 B(U)=e^{-aU^2/4+U/2+a\epsilon_0U}
        \sum_{d\mid D}e^{aq_d-\bar\sigma(U,a)\log d}.
\]
The independently checked inequality
$-aU_0/2+1/2+a\epsilon_0<0$ makes each forcing factor decrease for
$U\ge U_0$.  Thus, for $Z=\mathcal H+\mathcal R$,
\[
 Z'\le-\kappa Z+B(U_0)+C_0.                                  \tag{ET8}
\]

The primary $192$-bit Arb reconstruction sums all $2391233$ prefix
terms for both required logarithmic moments and all $64$ signed
divisor contributions.  It evaluates (ET7) by its exact antiderivative:
writing $B_0=\bar\sigma(U,a)-1$, it is
\[
 \frac{\sqrt\pi}{\sqrt a}e^{-B_0^2/a}
 \operatorname{erfi}\!\left(\frac{av/2-B_0}{\sqrt a}\right).
\]
The resulting fresh strict bounds are
\[
 Z(U_0)<0.981031775189216<\frac{491}{500},\qquad
 \frac{B(U_0)+C_0}{\kappa}<0.752526<\frac{491}{500}.
                                                               \tag{ET9}
\]

\subsection{An independent reconstruction of the complete Euler arithmetic}

A separate $256$-bit implementation groups by the represented integer
$m$, rather than by prefix sums, and evaluates
\[
 \beta_m(a)m^{-\bar\sigma(U_0,a)}
 =\sum_{d\mid\gcd(m,D)}\mu(d)
  \exp\!\left(a\left[q_d+\frac{(\log m-\log d)^2}{4}\right]
              -\bar\sigma(U_0,a)\log m\right)                 \tag{ET10}
\]
for every $2\le m\le K$.  It treats the six selected prime coefficients
as their exact algebraic zeros.  Every other one of the $2391226$
complete coefficients is independently enclosed strictly above zero.
Their full sum and the complete weighted correction are contained in
the primary enclosures.

The same implementation independently encloses each integral in
(ET7) by $512$ equal closed panels.  For
$g(v)=e^{av^2/4-(\bar\sigma-1)v}$ one has
\[
 g''(v)=\left[\frac a2+
        \left(\frac{av}{2}-(\bar\sigma-1)\right)^2\right]g(v)>0.
\]
The composite midpoint sum is a lower bound and the composite
trapezoid sum is an upper bound.  All $32768$ panels are evaluated with
directed Arb arithmetic.  They enclose the primary erfi value and give
the independent bound
$Z(U_0)<0.981031860662434<491/500$.
The forcing is separately recomputed as the positive six-factor
product
\[
 B(U_0)=e^{-aU_0^2/4+U_0/2+a\epsilon_0U_0}
       \prod_{p\in\mathcal P}
           \left(1+e^{a\log^2p/4-\bar\sigma(U_0,a)\log p}\right).
\]
Its independent equilibrium is less than $0.752526$.
These computations use different complete sum groupings and different
integral bounds, while sharing Arb/FLINT directed arithmetic.

Solving the differential inequality (ET8) gives
\[
 Z(U)\le e^{-\kappa(U-U_0)}Z(U_0)
       +(1-e^{-\kappa(U-U_0)})\frac{B(U_0)+C_0}{\kappa}
       <m_*:=\frac{491}{500}.                                \tag{ET11}
\]
By (ET2)--(ET7), this bounds the entire nonconstant coefficient mass
of $EA_N$ at the actual real part $\sigma$, uniformly for every
$a\le t\le b$ and every $N\ge H_*$.  This is a continuum proof on
the unbounded range.

\subsection{Factor bounds and the actual phase length}

At the adverse corner define
$q_p=e^{a\log^2p/4-\bar\sigma(U_0,a)\log p}$.
The fresh computation proves $0<q_p<1$ for all six primes.  Their
majorising amplitudes decrease with $U$ and $t$, so the corner bounds
are uniform. With $E=\rho e^{i\phi}$ they imply
\[
 \rho\ge\prod_p(1-q_p)>0.472061847298902>\frac{23}{50},\qquad
 \rho\le\prod_p(1+q_p)<1.82364<\frac{19}{10}.
                                                               \tag{ET12}
\]
In particular $E$ is nonzero and a differentiable phase $\phi$ may be
chosen locally; its derivative is unambiguous.

At $s=(1+ix)/2$, direct differentiation gives
\[
 \Re\alpha'(s)=\frac{7x^2-5}{(1+x^2)^2},\qquad
 \Im\alpha'(s)=-\frac{x(x^2+5)}{(1+x^2)^2},
\]
\[
 \sigma_x=\frac{tx(x^2+5)}{4(1+x^2)^2},\qquad
 \tau_x=\frac12+\frac{t(7x^2-5)}{4(1+x^2)^2}.                 \tag{ET13}
\]
For $Y=H_*^2$, $x>Y$ implies the exact bounds
\[
 0<\tau_x-1/2\le\frac{7b}{4Y^2}<e_0:=10^{-24},
\]
\[
 |w_x|\le\tau_x+\sigma_x
 <\frac12+\frac{7b}{4Y^2}+\frac b{4Y}+\frac b{Y^3}
 <c_0:=\frac{500000000001}{10^{12}}.
\]
Since $\phi_x=\Im(E_x/E)$,
\[
 |\phi_x|\le c_0\sum_p\frac{q_p\log p}{1-q_p},\qquad
 |E_x|\le c_0\prod_p(1+q_p)
                     \sum_p\frac{q_p\log p}{1+q_p}.
\]
The actual directed corner bounds give
\[
 L_1:=\sum_p\frac{q_p\log p}{1-q_p}
       <0.896716206645132<\frac9{10},\qquad
 |E_x|<0.579342126550494<\frac35.                              \tag{ET14}
\]

Let $\Theta$ be the real phase of $M_t(s)$, so $P=e^{i\Theta}$ and
$V=2\Re(PA_N)$. Since
$M_t(s)=M_0(s)e^{t\alpha(s)^2/4}$ and $M_0'/M_0=\alpha$, one has
\[
 4\Theta_x=2\Re\alpha+t\Re(\alpha\alpha')
           =4c(x)+t\Re(\alpha\alpha').
\]
The explicit formulas give
\[
 \Im\alpha=\tfrac12\arctan x-\frac{3x}{1+x^2}
        >\frac34-\frac7{2x}>0,
 \quad-\Im\alpha'>\frac1x,
 \quad\Re\alpha>0,\quad\Re\alpha'>0.
\]
Also $\log(1+r)\ge r/(1+r)$ for $r\ge0$ gives
\[
 4c(x)-\log(x/(4\pi))
 =\tfrac12\log(1+x^{-2})-\frac2{1+x^2}>-\frac3{2x^2}.
\]
It follows that
\[
 4\Theta_x-\log(x/(4\pi))
 >\frac{3a}{4x}-\frac{7b/2+3/2}{x^2}>0,                      \tag{ET15}
\]
where the final sufficient endpoint gate
$3aH_*^2/4>7b/2+3/2$ is checked exactly.

Set $\theta=\Theta-\phi$, $h=\theta_x$, and
$\ell_{\mathrm{actual}}=2h/\tau_x$.
Equations (ET2), (ET13)--(ET15) give
\[
 \ell_{\mathrm{actual}}
 >\ell_0(U,t):=
 \frac{2U-\epsilon_0-
       4c_0\sum_p q_p(U,t)\log p/(1-q_p(U,t))}{1+2e_0},       \tag{ET16}
\]
where $q_p(U,t)=e^{t\log^2p/4-\bar\sigma(U,t)\log p}$.
At the corner the fresh bound is
\[
 \ell_0(U_0,a)>27.588526177342325>27,
 \qquad\ell_0(U_0,a)>U_0+\log D.
\]
Each $q_p(U,t)$ decreases with both variables. Therefore
$\ell_0(U,t)-(U+\log D)$ increases with $t$ and has $U$-derivative
greater than $2/(1+2e_0)-1>0$.
Consequently $\ell_{\mathrm{actual}}>\log(DN)$ throughout the tail.
Since $\tau_x>1/2$, one also has $4h>\ell_{\mathrm{actual}}>27$.

\subsection{The completed-resultant matrix and every represented product}

Define
\[
 \mathcal T=\rho V/2-\frac{i}{2h}(\rho V)_x.
\]
As $\rho V=2\Re(e^{i\theta}EA_N)$, the constant term of $EA_N$
contributes exactly $e^{i\theta}$ to $\mathcal T$.
Consider any other represented term
$z=e^{i\theta}c_m(t)m^{-w}$, where $c_m(t)$ is independent of $x$
when $N$ and $t$ are fixed. Its derivative is
$z_x=(ih-w_x\log m)z$. Thus its two real coordinates in $\mathcal T$
are obtained from $(\Re z,\Im z)$ by
\[
 \begin{pmatrix}
 1&0\\2\varepsilon q&1-2q
 \end{pmatrix},\qquad
 \varepsilon=\frac{\sigma_x}{\tau_x},\qquad
 q=\frac{\log m}{\ell_{\mathrm{actual}}}\in[0,1].             \tag{ET17}
\]
The interval for $q$ includes every represented $m\le DN$, by (ET16),
without truncating the completed product back to $N$.

The matrix in (ET17) is a convex combination of the identity and
$S_\varepsilon=\left(\begin{smallmatrix}1&0\\2\varepsilon&-1\end{smallmatrix}\right)$.
The larger eigenvalue of $S_\varepsilon^TS_\varepsilon$ is
$(\sqrt{1+\varepsilon^2}+|\varepsilon|)^2$. Its operator norm is
therefore $\sqrt{1+\varepsilon^2}+|\varepsilon|$, and convexity of the
operator norm gives the same bound for (ET17).
Applying the triangle inequality to the complete nonconstant mass
in (ET11) proves
\[
 |\mathcal T-e^{i\theta}|
  <m_*(\sqrt{1+\varepsilon^2}+|\varepsilon|).                 \tag{ET18}
\]

\subsection{Physical real-jet conversion with exact rational reserve}

For the actual real normalisation $u=H_t/G$, define
\[
 \mathcal J=\rho u/2-\frac{i}{2h}(\rho_xu+\rho u_x).
\]
Since $u,V,u_x,V_x$ are real and $|\rho_x|\le|E_x|$,
(ET1), (ET12), (ET14) and $4h>27$ give
\[
 |\mathcal J-\mathcal T|^2\le
 D_2:=\left(\frac{19v}{20}\right)^2+
 \left[\frac2{27}\left(\frac{3v}{5}+\frac{19d}{10}\right)\right]^2
 =\frac{234549697}{116640000000000000000}.                     \tag{ET19}
\]
Equation (ET13) also gives
\[
 |\varepsilon|\le\varepsilon_*
 :=2b\left(\frac1{4Y}+\frac1{Y^3}\right),\qquad
 \iota:=m_*(\varepsilon_*+\varepsilon_*^2/2).
\]
The inequality $\sqrt{1+z^2}\le1+z^2/2$ bounds the excess norm in
(ET18) by $\iota$. Independent exact rational arithmetic proves
\[
 \iota<\frac1{500000},\qquad
 D_2<\left(\frac1{500000}-\iota\right)^2.
\]
Thus the sum of the physical error and completed-product norm inflation
is below $1/500000$, and
\[
 |\mathcal J|>1-\frac{491}{500}-\frac1{500000}
                    =\frac{8999}{500000}.                    \tag{ET20}
\]
Finally the same factor bounds imply
\[
 |\mathcal J|
 \le\left(\frac{19}{20}+\frac2{27}\frac35\right)|u|
       +\frac2{27}\frac{19}{10}|u_x|
 =\frac{179}{180}|u|+\frac{19}{135}|u_x|.
\]
Cauchy--Schwarz therefore gives the strict physical bound
\[
 u^2+u_x^2>
 \frac{(8999/500000)^2}{(179/180)^2+(19/135)^2}
 =\frac{59035878729}{183840625000000}>\frac3{10000}.            \tag{ET21}
\]
Since $G(x)>0$, $H_t(x)=H'_t(x)=0$ would force $u(x)=u_x(x)=0$,
contradicting (ET21).

\subsection{The closed outer seam}

For $x_j(t)=4\pi(j^2-t/16)$ the independent $256$-bit seam check proves
\[
 W-x_{H_*}(a)>0.2869527460786004>0,
 \qquad x_{H_*+1}(b)-W>60318591.19442919>0.
\]
It separately proves $W-1<4\pi H_*^2<W$, verifying the exact ceiling.
Hence $N(t,W)=H_*$ for every $a\le t\le b$, and every $x\ge W$
has natural index at least $H_*$. All bounds above consequently hold
on the closed physical half-line $x\ge W$.
The independently certified finite exterior includes $X\le x\le W$,
so the two domains overlap at $x=W$ with no omitted natural-index cell.
Evenness gives the corresponding negative exterior $x\le-W$.

\section{Critical-time compactness and proof of the main theorem}
\label{sec:descent}

\begin{lemma}[Uniform remote reality on a positive time band]
\label{lem:uniform-remote-band}
For every $0<r\le b\le1/2$ there is a finite $R_0=R_0(r)$ such that
\[
 t\in[r,b],\quad H_t(z)=0,\quad |\Re z|\ge R_0
 \quad\Longrightarrow\quad z\in\mathbb R.
\]
\end{lemma}
\begin{proof}
Polymath's Theorem~1.5(i) \cite{Polymath2019} supplies an absolute
constant $C>0$ for the positive-real-part cutoff $e^{C/t}$ throughout
$0<t\le1/2$. Choose $R_0=e^{C/r}$. Then $R_0\ge e^{C/t}$ for
$r\le t\le b$, so the cited result applies. Evenness handles negative
real parts.
\end{proof}

The assertion is uniform only after a positive lower time $r$ is fixed.
No uniform cutoff as $t\downarrow0$ is asserted. The proof uses the
reality conclusion of the published theorem; its asymptotic location
formula and its separate uniqueness statement are not needed.

\begin{lemma}[Finite collision at a positive critical time]\label{lem:stall}
If $0<\Lambda\le1/2$, then $H_\Lambda$ has a multiple real zero at a
finite point.
\end{lemma}
\begin{proof}
By the defining de Bruijn--Newman theorem, every zero of $H_\Lambda$ is
real. Suppose, for a contradiction, that every such zero is simple.
Apply Lemma~\ref{lem:uniform-remote-band} to
$[r,\Lambda]$, where $r=\Lambda/2>0$. Choose $R>R_0$ with
$H_\Lambda(R)\ne0$ and $H_\Lambda(-R)\ne0$. Such a choice is possible
because $H_\Lambda$ is a nonzero entire function and its zeros are
discrete. Nontriviality also follows directly from
$H_\Lambda(0)>0$.

Let $K=\{z:|\Re z|\le R,\ |\Im z|\le2\}$. Its boundary is zero-free
at time $\Lambda$: the vertical sides can meet the real axis only at
$\pm R$, and all zeros at this time are real. It contains finitely many
zeros. Put pairwise disjoint closed disks around them, centred on the
real axis and contained in the interior of $K$, with exactly one zero
in each disk. Their boundaries and the closed complement of their
interiors in $K$ are zero-free at $\Lambda$.

Joint continuity gives a zero-free complement for every real $t$
sufficiently close to $\Lambda$. Rouch\'e's theorem on the disk
boundaries preserves one zero per disk, counted with multiplicity.
Each disk is conjugation-invariant and $H_t$ is real-type, so its unique
zero is real. Thus every zero inside $K$ is real at every such time.

Choose one of these times with $r\le t<\Lambda$. The coarse de Bruijn
strip excludes zeros with $|\Im z|>1$. Every zero outside $K$ must
therefore satisfy $|\Re z|>R>R_0$ and is real by
Lemma~\ref{lem:uniform-remote-band}. This would make all zeros of $H_t$
real at a time below $\Lambda$, a contradiction. At least one zero of
$H_\Lambda$ is therefore multiple. It is real and finite.
\end{proof}

No assertion about remote simplicity is used in this proof. No numerical
value of $R_0$ is required, and no lower-bound theorem for $\Lambda$ is
used: positivity of $\Lambda$ will follow from the contradictory case
in the next theorem.

\begin{theorem}[Regional assembly criterion]\label{thm:regional-assembly}
Set
\[
 a=\frac{156821}{1000000},\qquad b=\frac15,\qquad
 X=6000000185827,\qquad W=72382294738709.
\]
Assume the following three physical statements:
\begin{enumerate}
\item \textbf{Closed stationary boundary.}
      $H_t(X+iy)\ne0$ for every $(t,y)\in[0,b]\times[0,1]$.
\item \textbf{Closed finite exterior.}
      $(H_t(x),\partial_xH_t(x))\ne(0,0)$ whenever
      $a\le t\le b$ and $X\le x\le W$.
\item \textbf{Closed infinite exterior.}
      $(H_t(x),\partial_xH_t(x))\ne(0,0)$ whenever
      $a\le t\le b$ and $x\ge W$.
\end{enumerate}
Then, using the adopted published Platt--Trudgian results and the
Polymath strip and remote-reality theorems,
\[
 \Lambda<\frac{156821}{1000000}.
\]
\end{theorem}

\begin{proof}
Platt--Trudgian's Corollary~2 \cite{PlattTrudgian2021} gives
$\Lambda\le b$. If $\Lambda<a$, the conclusion already holds.
Otherwise $a\le\Lambda\le b$, so $0<\Lambda\le1/2$ and
Lemma~\ref{lem:stall} supplies a multiple real zero of $H_\Lambda$.

The first hypothesis and Corollary~\ref{cor:physical-interior} exclude
all multiple real zeros on $|x|\le X$ for every $0<t\le b$.
The second and third hypotheses exclude them on $[X,W]$ and
$[W,\infty)$ for $a\le t\le b$. Evenness gives the same exclusion
at negative $x$, since $H_t(-x)=H_t(x)$ and
$\partial_xH_t(-x)=-\partial_xH_t(x)$.
These closed regions cover the entire real line. Hence the multiple
zero at $t=\Lambda\in[a,b]$ cannot exist, a contradiction.
\end{proof}

\paragraph{What the criterion requires from the certificates.}
The finite-exterior hypothesis is a statement about the physical real
jet. A lower bound for a complex source polynomial alone does not imply
it. Its proof must include the carrier derivative, both remainder
estimates, all moving-index and source-handoff terms, and a strictly
positive final Hermitian reserve on the stated closed domain.
Similarly, the separator must be transferred to $H_t$ with its full
complex error allowance, and the infinite-exterior bound must imply
nonvanishing of the physical pair $(H_t,\partial_xH_t)$.

If a regional calculation gives $u^2+u_x^2>0$ for
$u(t,x)=H_t(x)/A(t,x)$ with a nonzero real differentiable factor $A$,
then it does imply the required physical exclusion: at any common zero
of $H_t$ and $\partial_xH_t$, both $u$ and $u_x$ vanish. The factors used
in different regions need not coincide, because the hypotheses concern
the same physical function after this transfer.

Both time endpoints and both spatial seams are included in the criterion.
It excludes $\Lambda=a$ as well as $a<\Lambda\le b$, which is why the
conclusion is strict. This is a qualitative strict inequality; no
specific further decrement below $a$ is furnished by compactness.

The preceding theorem isolates the analytic assembly from the regional
calculations. Its hypotheses are discharged by the specific modules of
this manuscript in the following application.

\begin{theorem}[Computer-assisted upper bound]
\label{thm:global}
In the normalisation $H_0(z)=\xi((1+\ii z)/2)/8$, the exact-input
certificates and the analytic transfers developed in this manuscript give
\[
 \boxed{\Lambda<\frac{156821}{1000000}=0.156821.}
\]
The adopted published inputs are precisely those listed in
Section~\ref{sec:normalisation}; their original large computations are
not rerun here.
\end{theorem}
\begin{proof}
Section~\ref{sec:separator} proves the closed stationary boundary with
normalised physical modulus greater than $7/800$.
Theorem~\ref{thm:multi-gluing} and the five complete source certificates
give a stable full prime polynomial with modulus greater than $1/790$
on the closed finite exterior. Its coefficients agree exactly with the
balanced gauge through $c_{\rm pol}=\Rea\alpha/2$.
Corollary~\ref{cor:admission} therefore supplies the second physical
hypothesis, including every natural-index wall and every source handoff.
Section~\ref{sec:tail} proves $u^2+u_x^2>3/10000$ on the whole infinite
exterior with the same real gauge; this supplies the third hypothesis.
The three domains have the common time band $[a,b]$ where required,
and cover both spatial seams. Theorem~\ref{thm:regional-assembly} applies.
\end{proof}

The theorem is the mathematical claim of proposed version 0.3. The fresh
finite replay and the proof audit are documented separately; author approval
and independent human mathematical examination remain part of preparing a
submission. The statement does not represent completed peer review.

\section{Discussion and outlook}\label{sec:discussion}
The proof gives a strict upper bound from a closed-band exclusion of
real multiple zeros rather than a canopy of complex nonvanishing.
This is important independently of the decimal endpoint: it separates
the global de Bruijn--Newman argument into a compact initial-region
theorem, finite spatial stability certificates, and an infinite-tail
physical differential inequality. The same published finite-height
RH calculation is sufficient for the present approach, and no
independent hypothesis about the location of unverified zeta zeros
is made.

The four-prime Cauchy--Schwarz obstruction identifies a limit of the
particular positive-weight Gram source family, not a lower limit on
$\Lambda$. Enlarging the smooth-prime basis or using a different stable
polynomial architecture may alter that feasibility constraint.
The arithmetic and interval-computation requirements of such an
extension remain substantial, and the present paper makes no further
numerical claim. The full certificate record is important because the
strict inequality relies on outward-rounded bounds at finite,
closed geometries rather than on sample values of the model function.

\section{Computational evidence and independent reproducibility}\label{sec:replay}
The proof uses a finite collection of inequalities for explicitly defined
Dirichlet-polynomial coefficients and parameter rectangles. The
computation is designed to be \emph{proof-bearing}: each finite box is
closed, interval operations are outward rounded, and every interval
lower bound includes approximation and discarded-support charges.
The analytic argument invokes the certified statements themselves,
not an informal plot or a sampled numerical minimum.

In the frozen computational reconstruction supplied by Raymond, the
unchanged reference driver \texttt{replay\_main.py} reconstructs five
full Fourier source tensors (664,505 structural coefficients), checks
44,182 torus boxes, rebuilds all 184 Acb separator moments and
17,385 separator rectangles, and evaluates the 56 recovery inequalities
with separate geometry checks. The accompanying record contains source
hashes, full execution logs, exact-rational gauge tests, centred-moment
checks and complete Euler-tail arithmetic. Independent diagnostic
programs check portions of the same chain by alternative formulas;
most directed computations share the Arb/FLINT numerical library.
Accordingly, those diagnostics are not described as an independently
implemented, end-to-end interval proof.

The complete computational certificates, reconstruction programs, execution
logs, manifest files and replay instructions are publicly available in the
versioned Planat--Raymond reproducibility repository
\cite{PlanatRaymondCertificates2026}:
\begin{flushleft}\small\url{https://github.com/Dale-Dug/planat-raymond-collision-exclusion-0156821}\end{flushleft}
The fixed release is
\begin{flushleft}\small\url{https://github.com/Dale-Dug/planat-raymond-collision-exclusion-0156821/releases/tag/certificates-0156821-v1}\end{flushleft}
under tag \texttt{certificates-0156821-v1}, at source commit
\href{https://github.com/Dale-Dug/planat-raymond-collision-exclusion-0156821/commit/18dbc01ea36bf9a24d8f79822b5ce06b545e4ace}{\texttt{18dbc01ea36b}} (full identifier in the ancillary README).
The complete archive is
\path{Planat_0156821_Reproducibility_20261104.zip}
with SHA-256 checksum
\begin{quote}
{\footnotesize\ttfamily 08bf5123ac79a438b626275eb8174df6d8cdfd069e7c7a0a70c1c4807b1ebde7}.
\end{quote}
The archived release preserves the frozen v0.3 scientific inputs and
manuscript v0.6 for provenance; this revision changes the repository
citation and the AI-assistance declaration, not any scientific input or
certificate. The present arXiv source contains compact ancillary
proof-input hashes and representative checking scripts, while the
complete computational archive and its reproduction instructions are
available through the release. A file-integrity check is not a
substitute for re-executing the directed-interval computations.

\paragraph{Separation from exploratory extensions.}
The structural obstruction for the four-prime Gram family explains why
the chosen endpoint lies close to that family's feasibility limit.
It does not prohibit a different prime basis or a better lower bound.
Preliminary five-prime candidates toward $0.155$ are outside the
present theorem and no assertion about them enters the proof.

\section*{Declaration of AI assistance}
AI-based systems were used in mathematical exploration, program construction,
checking and exposition. We used OpenAI ChatGPT 6.1. An adversarial audit was
conducted using Anthropic Claude Opus 5.5.

\end{document}